\documentclass[12pt]{article}
\usepackage{amsmath}
\usepackage{amssymb}
\usepackage{amsthm}
\usepackage{mathrsfs}
\usepackage{authblk}

\newcommand{\R}{\mathbb R}

\newcommand{\C}{\mathbb C}
\newcommand{\beq}{\begin{eqnarray}}
\newcommand{\eeq}{\end{eqnarray}}
\newcommand{\beqst}{\begin{eqnarray*}}
\newcommand{\eeqst}{\end{eqnarray*}}

\title{\bf Uniqueness for an inverse problem\\ for a semilinear time-fractional\\ diffusion equation}
\author{Jaan Janno\footnote{Corresponding author}, Kairi Kasemets}
\affil{\small Tallinn University of Technology  Ehitajate tee 5,\break 19086 Tallinn, Estonia\break\break E-mails: J.Janno: jaan.janno@ttu.ee; K. Kasemets: kairi.kasemets@ttu.ee}

\date{}

\begin{document}

\maketitle

\abstract{\small An inverse problem to determine a space-dependent factor in a semilinear time-fractional diffusion equation is considered.
Additional data are given in the form of an integral with the Borel measure over the time.
Uniqueness of the solution of the inverse problem is studied.
The method uses a positivity principle of the corresponding differential equation that is also proved in the paper.}

\section*{Introduction}

Anomalous diffusion processes in porous, fractal, biological etc. media are described by differential equations containing fractional derivatives. Depending on the nature of the process, the model may involve either fractional time or fractional space derivatives or both ones \cite{Gor,Mag,Zas}.

In many practical situations the properties of the medium or sources are unknown and they have to be reconstructed
solving inverse problems \cite{Che,Jin,SakYam,Zhang}. Then some additional information
  for the solution of the differential equation is needed to recover
the unknowns.
If space-dependent quantities are to be determined, such  additional conditions may involve instant (e.g. final)
 measurements or integrated measurements over time.

In case space and time variables are separable, the solution of the time-fractional diffusion equation can be expressed by a formula that
is deduced by means of the Fourier expansion with respect to the space variables and integration using Mittag-Leffler functions
with respect to the time variable. This formula enables to prove uniqueness of reconstruction of space-dependent factors
 of source terms from final overdetermination \cite{Fu,Kir,Wei}. However, this method fails even in the linear case
  when the variables are not separable.

In the present paper we prove the uniqueness for an inverse problem to determine a space-dependent factor in
 the time-fractional diffusion equation in a more general
case when the equation may be semilinear. Additional condition is given in a form of an integral with  Borel measure over the time that includes as a
particular case the final overdetermination. Such an inverse problem has possible applications in modelling of
fractional reaction-diffusion processes \cite{Gaf,Kaz,Rida}, more precisely in reconstruction of certain parameters of
inhomogeneous media.

Our results are global in time, but  contain certain cone-type restrictions that may depend on a time interval.
We will adjust a method that was
  applied to inverse problems for usual parabolic equations \cite{Is1,Is2} and generalized
   to an integrodifferential case  \cite{JaKa} and parabolic semilinear case \cite{Ber}.
    The method is based on positivity principles of  solutions of differential equations that
follow from extremum principles. Extremum principles for time-fractional diffusion equations
 are proved in the linear case in \cite{Bru,Lu1,ReLu} and in a nonlinear divergence-type case in \cite{Zach1}.
  But these results are not directly applicable
in our case, because of the lack of the semilinear term. Therefore, we prove an independent positivity principle in our paper.

To the authors' opinion, such a positivity principle has a scientific value independently of the inverse problem, too.
States of several reaction-diffusion models are  positive functions, e.g.  probability densities
\cite{Kaz}.

The stability of the solution of the inverse problem will not be studied in this paper.
In the linear case the stability follows from the uniqueness by means of the Fredholm alternative \cite{SakYam}. The solution
continuously depends on certain derivatives of the data, hence the problem is moderately ill-posed.

The plan for the paper is as follows. In the first section we formulate the direct and inverse problems. Second section contains
auxiliary results about a linear direct problem that are necessary for the analysis of the inverse problem. Third section is
devoted to the positivity principle. In the fourth section we prove the main uniqueness theorem in case of the general additional condition
involving the Borel measure. The next section contains a particular uniqueness result in a case of the Lebesgue measure with a weight. In the last
section we discuss some crucial assumptions of the uniqueness theorems.

\section{Problem formulation}\label{s:1}

Let $\Omega\subset\R^n$ be a bounded open domain with a sufficiently
smooth boundary $\partial \Omega$. We
consider the semilinear fractional diffusion equation
\beq\label{maineq} D_t^\beta [u(t,x)-u_0(x)]\, =\, A(x)u(t,x) +
f(u(t,x),t,x)\, ,\quad t\in (0,T),\; x\in \Omega\, , \eeq
with the initial  and boundary conditions
\beq\label{ini} &&u(0,x)=u_0(x)\, ,\quad x\in\Omega,
\\[1ex] \label{bound}
&&{\cal B}u(t,x)=g(t,x)\, ,\;\; (t,x)\in (0,T)\times {\partial\Omega}
\eeq
where
\beq\label{Bdef}
\mbox{either\,\, (I)}\;\;{\cal B}u=u\quad\mbox{or\,\,\,\, (II)}\;\;{\cal B}u=\omega\cdot \nabla u
\eeq
with some $x$-dependent vector function $\omega(x)=(\omega_1(x),\ldots,\omega_n(x))$ such that $\omega(x)\cdot \nu(x)>0$
where $\nu$ is  the outer normal
to ${\partial\Omega}$. Here and in the sequel $\nabla$ denotes the gradient operator with respect to the space variables.
Moreover,
$$D_t^\beta v(t)\,=\, {1\over\Gamma(1-\beta)}{d\over dt}\int_0^t (t-\tau)^{-\beta}
v(\tau)d\tau$$
is the Riemann-Liouville fractional derivative of order $\beta\in (0,1)$
and the operator $A$ has the form
 $$
 A(x) = \sum_{i,j=1}^n a_{ij}(x){\partial^2\over\partial x_i\partial x_j}+\sum_{j=1}^n a_{j}(x){\partial\over\partial x_j},
 $$
where the principal part is uniformly elliptic, i.e.
$$
\sum_{i,j=1}^n a_{ij}(x)\xi_i\xi_j\ge c|\xi|^2\quad \forall \xi\in\R^n,\, x\in\Omega\quad\mbox{for some $c>0$.}
$$

Note that in case of sufficiently smooth $u$ the term in the left-hand side of \eqref{maineq} is actually the Caputo derivative of $u$, i.e.
$D_t^\beta [u(t,x)-u_0(x)]={1\over\Gamma(1-\beta)}\int_0^t (t-\tau)^{-\beta}u_\tau (\tau,x)d\tau$.
The conditions \eqref{maineq}, \eqref{ini} and \eqref{bound} form a direct problem for the function $u$.

The existence and uniqueness of the solution of
\eqref{maineq} - \eqref{bound} in case
$f(w,t,x)$ depends only on $w$ are proved in \cite{Kras,Lu2}. More precisely, \cite{Lu2} considers the problem with Dirichlet boundary conditions
and \cite{Kras} the one-dimensional problem with Neumann boundary conditions.
In the appendix of the paper we prove independent existence and uniqueness theorems for
\eqref{maineq} - \eqref{bound} in the general form.

\smallskip
We will consider the case when the nonlinearity function $f$ has the following form:
\beq\label{fform}
f(w,t,x)=a(w,t,x)z(x)+b(w,t,x),
\eeq
where $a$ and $b$ are given but the factor $z$ is unknown.

Let us formulate an inverse problem.
\\[1ex]
{IP}.\, Let $\mu$ be a positive Borel measure such that ${\rm supp}(\mu)\cap (0,T]\ne \emptyset$.
Determine a pair of functions $(z,u)$  such that the conditions  \eqref{maineq}, \eqref{ini}, \eqref{bound}, \eqref{fform} and
 the  additional condition
\beq\label{add}
\int_0^T u(t,x)d\mu\, =\, d(x)\, ,\quad x\in\Omega
\eeq
with some given function $d$ is satisfied.

\smallskip
We note that  a special case of $\mu$ is the Dirac measure concentrated at the final moment $t=T$. Then the condition \eqref{add} reads
$u(T,x)=d(x)$, $x\in\Omega$.

\section{Preliminaries}\label{s:2}
\setcounter{equation}{0}

In this section we formulate and prove some auxiliary results. In addition to $D_t^\beta$,
we introduce the
  operator of fractional integration of order $\gamma>0$:
$$
J_t^\gamma w(t)={1\over \Gamma(\gamma)}\int_0^t (t-\tau)^{\gamma-1}w(\tau)d\tau.
$$

In the sequel we consider the fractional differentation and integration in Bessel potential and H\"older spaces.
We consider the abstract Bessel potential spaces $H_p^\beta([0,T];X)$ for Banach spaces  $X$  of the class
${\cal HT}$ that is defined in the following manner
 (\cite{zachdiss}, p. 18,  \cite{pruess}, p. 216):\footnote{For example, $L_s(\Omega)$, $1<s<\infty$, are of the class
${\cal HT}$.}
$$\begin{array}{ll}
&{\cal HT}=\{X\; :\; X - \mbox{Banach space, the Hilbert transform is bounded in $L^q(\R;X)$}
\\
&\qquad \mbox{for some $q\in (1,\infty)$}\}
\end{array}
$$

Due to an embedding theorem, it holds $H_p^\beta([0,T];X)\subset C^\alpha([0,T];X)$ provided
${1\over\beta}<p<\infty$, $\alpha\in (0,\beta-{1\over p})$ and $X$ is of the class
${\cal HT}$.

Let us formulate four lemmas that directly follow from known results in the literature.
\\[1ex]
\noindent{\bf Lemma 1}. {\it Let $X$ be a Banach space of the class
${\cal HT}$ and ${1\over\beta}<p<\infty$.
If $w\in H_p^\beta([0,T];X)$ and $w(0)=0$ then $D_t^\beta w\in L_p((0,T);X)$.
If $w\in L_p((0,T);X)$ then $J_t^\beta w \in H_p^\beta([0,T];X)$ and $J_t^\beta w(0)=0$.
Moreover,
 the relations
\beq\label{JD}
&&J_t^\beta D_t^\beta (w-w(0))\, =\, w-w(0)\quad\forall w\in H_p^\beta([0,T];X),
\\
\label{DJ}
&&D_t^\beta J_t^\beta w\, =\, w \quad\forall w\in L_p((0,T);X)
\eeq
are valid.}
\\[1ex]
{\it Proof}. These assertions follow from arguments presented in \cite{zachdiss}, p. 28, 29.
\\[1ex]
{\bf Lemma 2}.
{\it
Let $X$ be a Banach space. If $w\in C^\alpha ([0,T]; X)$, $\alpha\in (0,1)$, $w(0)=0$
and $\beta\in (0,1-\alpha)$ then $J_t^\beta w\in C^{\alpha +\beta}([0,T]; X)$ and $J_t^\beta w(0)=0$.
If $w\in C^\alpha ([0,T]; X)$, $\alpha\in (0,1)$, $w(0)=0$ and $\beta\in (0,\alpha)$ then
$D_t^\beta w\in C^{\alpha -\beta}([0,T];X)$ and $D_t^\beta w(0)=0$.
}\\[1ex]
\noindent
{\it Proof.} Assertions of the lemma  follow from  Thms 14, 19 and the first
part of Thm. 20 of \cite{HaLi} if we continue $w(t)$ by zero for $t<0$.
Although \cite{HaLi} considers the case $X=\R$, all arguments included in  proofs of the mentioned
theorems automatically hold in case of arbitrary Banach space $X$, too.
\hfill $\Box$
\\[1ex]
{\bf Lemma 3}
{\it Let ${1\over\beta}<p<\infty$, $p\not\in\{{1\over\beta}+{1\over 2};{2\over \beta}+1\}$.
Assume $a_{ij}\in C(\overline\Omega)$, $a_j\in L_\infty(\Omega)$, $q\in L_\infty((0,T)\times \Omega)$,
$\omega\in (C^{r}({\partial\Omega}))^n$, $r>1-{1\over p}$.
Then the problem
\beq\label{l2prob}
&&\hskip -1truecm w(t,x)=J_t^\beta (A(x)+q(t,x))w(t,x)+\varphi(t,x)\, ,\quad x\in \Omega\, ,\; t\in (0,T),
\\ \label{l2bound}
&&{\cal B}w(t,x)=0\, ,\;\; (t,x)\in
(0,T)\times {\partial\Omega}
\eeq
has unique solution
\beq\label{Updefi}
w\in U_p:=H_p^\beta([0,T];L_p(\Omega))\cap L_p((0,T);H_p^2(\Omega))
\eeq
such that $w(0,\cdot)=0$ if and only if $\varphi\in \Phi_p=\{\varphi\in H_p^\beta([0,T];L_p(\Omega))\, :\, \varphi(0,\cdot)=0\}$. The
operator $\mathcal S$ that maps  $\varphi$ to $w$ is continuous from the space $\Phi_p$
to the space $U_p$.}
\\[1ex]
{\it Proof.} The existence and uniqueness assertions of Lemma 3 are a particular case of a
more general maximal regularity result proved in Thm 4.3.1 of \cite{zachdiss}. The proof of the cited theorem
is based on approximation of the problem by a sequence of localized problems with constant coefficients.
For the localized problems, existence and uniqueness theorem for an abstract parabolic evolutionary integral equation containing
unbounded operator is applied.
The study of the such an abstract equation is based on the construction of a solution in the form of a
variation of parameters formula that contains a convolution of a resolvent (operator) of the equation with the given right-hand side.
The desired assertions follow from the properties of the resolvent. A formula of the resolvent is constructed
and its properties established by means of the Laplace transform. The analysis of the mentioned abstract parabolic equation
is contained also in \cite{Zach2}.
 \hfill $\Box$
\\[1ex]
{\bf Lemma 4} {\it Let $X$ be a Banach space and ${\cal A}\, :\, D({\cal A})\subset X\to X$ a closed
densely defined unbounded operator satisfying the
following property:
\beq\label{sect}
\rho({\cal A})\supset \Sigma(\beta\pi/2)\, ,\;\;
\exists M>0\, :\, \|(\lambda-{\cal A})^{-1}\|\le {M\over |\lambda|}\;\; \forall \lambda\in \Sigma(\beta\pi/2)
\eeq
where $\rho({\cal A})$ is the resolvent set of ${\cal A}$ and
$\Sigma(\theta)=\{\lambda\in\C\, :\, |{\rm arg}\lambda|<\theta\}$. Let $X_{\cal A}$ be the domain of ${\cal A}$ endowed with the
graph norm $\|z\|_{X_{\cal A}}=\|z\|+\|{\cal A}z\|$.
If $\varphi\in C^\alpha ([0,T];X)$ for some $\alpha\in (0,1)$ and $\varphi(0)=0$ then the equation
$$
w(t)={\cal A}J_t^\beta w(t)+\varphi (t)\, ,\quad t\in [0,T]
$$
has a  solution $w\in C^\alpha ([0,T];X)$ satisfying the properties $w(0)=0$, $J_t^\beta w\in \break C^\alpha ([0,T];X_{\cal A})$.
The operator  ${\cal Q}_\alpha$ that maps  $\varphi$ to $w$ is continuous from
 the space $\{\varphi\, :\, \varphi\in C^\alpha ([0,T];X),\, \varphi(0)=0\}$
to the space $C^\alpha ([0,T];X)$.
If, additionally, $\varphi\in C^\alpha ([0,T];X_{\cal A})$ then $w\in C^\alpha ([0,T];X_{\cal A})$.}
\\[1ex]
{\it Proof}. Lemma 4 follows from Thm. 2.4, the estimate (2.31) in the proof of Thm. 2.4 and Example 2.1 of \cite{pruess}.
 The idea of the proof is analogous to the
proof of previous lemma. The solution is constructed in a form of a
convolution that contains a resolvent operator and properties of the resolvent are established by means of the Laplace transform.
\hfill $\Box$
\\[2ex]
\indent Now we  prove an  existence theorem in H\"older spaces for the direct problem \eqref{maineq}, \eqref{ini}, \eqref{bound}
in the linear case. This result  will be used in the analysis of IP.
\\[1ex]
{\bf Theorem 1}. {\it
Let $f(w,t,x)=q(t,x)w+\psi(t,x)$, $u_0=0$, $g=0$.
Assume that $a_{ij},a_j\in C^\alpha(\overline\Omega)$, $\omega\in (C^{1+\alpha}({\partial\Omega}))^n$  and
\beqst
&&\psi\in C^{\alpha+\beta}([0,T];L_p(\Omega))\cap C^{\alpha}([0,T];C^\alpha(\overline\Omega))\, ,\;\;\psi(0,\cdot)=0\, ,\\
&&q\in C^{\alpha+\beta}([0,T];L_\infty(\Omega))\cap C^{\alpha}([0,T];C^\alpha(\overline\Omega))
\eeqst
with some
$\alpha\in (0,1)$ and $p\in \bigl(\max\{{n\over 2}; {1\over\beta}\},\infty\bigl)$, $p\not\in\{{1\over\beta}+{1\over 2};{2\over \beta}+1\}$.
Then there exists $\alpha_1\in (0,1-\beta)$ such that the problem \eqref{maineq}, \eqref{ini}, \eqref{bound}
has a  solution
\beq\label{U}
u\in   C^{\alpha_1+\beta}([0,T];C^{\alpha_1}(\overline\Omega))\cap C^{\alpha_1}([0,T];C^{2+\alpha_1}(\overline\Omega))
\eeq
with $D_t^\beta u\in C^{\alpha_1}([0,T];C^{\alpha_1}(\overline\Omega))$.
The solution is unique in the wider space
$U_p$.
If in addition,
\beqst
\psi_t\in L_p((0,T);L_p(\Omega)),\; q_t\in L_p((0,T);L_p(\Omega))\cap L_1((0,T);L_\infty(\Omega))
\eeqst
then $u_t\in U_p$ and $u_t(0,\cdot)=0$.
}
\\[1ex]
{\it Proof.} Note that the problem \eqref{maineq}, \eqref{ini}, \eqref{bound}  is in the space
$U_p$
equivalent to the following problem:
\beq\label{tep1}
&&\hskip -1truecm u(t,x)=J_t^\beta (A(x)+q(t,x))u(t,x)+J^\beta_t\psi(t,x)\, ,\quad x\in \Omega\, ,\; t\in (0,T),
\\ \label{tep2}
&&{\cal B}u(t,x)=0\, ,\;\; (t,x)\in
(0,T)\times {\partial\Omega}\, .
\eeq
By Lemma 3,  the latter problem has a unique solution
$u\in U_p$. This proves the existence and uniqueness assertions in $U_p$.

Next we are going to prove the inclusion \eqref{U}. According to arguments of Section 3.1.1 and Thm 3.1.3 of \cite{Luna},
there exists $\xi\in\R$ such that the operator ${\cal A}=A+\xi$  with the domain $D({\cal A})=\{z\in W_p^2(\Omega)\, :\,
{\cal B}z|_{{\partial\Omega}}=0\}$ in $X=L_p(\Omega)$ satisfies the following conditions:
$$
\rho({\cal A})\supset \Sigma(\pi/2)\, ,\;\;
\exists M>0\, :\, \|(\lambda-{\cal A})^{-1}\|\le {M\over |\lambda|}\;\; \forall \lambda\in \Sigma(\pi/2).$$
Since $0<\beta<1$, this relation implies \eqref{sect}. Moreover, ${\cal A}$ is closed and
densely defined in $X$. Thus, ${\cal A}$ satisfies the assumptions of Lemma 4.
Let us consider the  problem
\beq\label{tep3}
&&\hskip -1truecm
v(t,x)= (A(x)+\xi)J_t^\beta v(t,x)+\varphi(t,x)\, ,
\quad x\in \Omega\, ,\; t\in (0,T),
\\ \label{tep4}
&&{\cal B}v(t,x)=0\, ,\;\; (t,x)\in
(0,T)\times {\partial\Omega}
\eeq
with $\varphi(t,x)=(q(t,x)-\xi)u(t,x)+\psi(t,x)$.
Observing that $u\in U_p$ and $U_p\subset C^{\alpha_2}([0,T];L_p(\Omega))$ with $\alpha_2\in (0,\beta-{1\over p})$,
by the embedding theorem, as well as the assumptions of the theorem we have
$\varphi\in C^{\alpha_2}([0,T];L_p(\Omega))$ and $\varphi(0,x)=0$.
Further, choosing some number $\alpha_3\in (0,\alpha_2]\cap (0,1-\beta)$ we also have
$\varphi\in C^{\alpha_3}([0,T];L_p(\Omega))$.
Applying Lemma 4 to the problem \eqref{tep3}, \eqref{tep4}, we conclude that it has a solution
$v\in C^{\alpha_3}([0,T];L_p(\Omega))$. Taking the operator $J_t^\beta$ from the relations
\eqref{tep3}, \eqref{tep4} and subtracting from \eqref{tep1}, \eqref{tep2} we see that the function $w=u-J_t^\beta v$ solves the problem
\eqref{l2prob}, \eqref{l2bound} with $\varphi=0$. Due to the uniqueness  it holds $w=0$, hence
$u=J_t^\beta v$. By the proved relation $v\in C^{\alpha_3}([0,T];L_p(\Omega))$
and Lemma 2 we get $u\in C^{\alpha_3+\beta}([0,T];L_p(\Omega))$.

According to the latter relation and assumptions of the theorem we can improve
properties of $\varphi$. Namely, it holds $\varphi\in C^{\alpha_4+\beta}([0,T];L_p(\Omega))$ with $\alpha_4=\min\{\alpha_3;\alpha\}\in (0,1-\beta)$.
Applying Lemma 4 again to \eqref{tep3}, \eqref{tep4}, we obtain
 $v\in C^{\alpha_4+\beta}([0,T];L_p(\Omega))$ with $J_t^\beta v\in C^{\alpha_4+\beta}([0,T];W_p^2(\Omega))$.
From the latter inclusion and the
embedding theorem we deduce
\beq\label{vahr}
u=J_t^\beta v\in C^{\alpha_4+\beta}([0,T];C^\gamma(\overline\Omega))\eeq with some $\gamma>0$. Lemma 2 implies
$D_t^\beta u\in C^{\alpha_4}([0,T];C^\gamma(\overline\Omega))$.

Further, let us rearrange the terms
in \eqref{maineq}, \eqref{ini}, \eqref{bound} to get the following family of elliptic problems:
\beq\label{tep7}
&&\hskip -1truecm A(x)u(t,x)=\phi(t,x)\, ,\quad x\in \Omega\, ,\; t\in (0,T),
\\ \label{tep8}
&&{\cal B}u(t,x)=0\, ,\;\; x\in
{\partial\Omega}\, ,\; t\in (0,T)
\eeq
with $\phi(t,x)=D_t^\beta u(t,x)-q(t,x)u(t,x)-\psi(t,x)$. According to the proved properties of $u$ and the assumptions of theorem
we have $\phi\in C^{\alpha_4}([0,T];C^{\alpha_5}(\overline\Omega))$ with $\alpha_5=\min\{\gamma;\alpha\}$. Using the well-known Schauder-estimates
for the
elliptic problems (e.g. \cite{Ladyz}, Ch. III) and taking the H\"older-estimates with respect to $t$ we obtain
$u\in C^{\alpha_4}([0,T];C^{2+\alpha_5} (\overline\Omega))$.
Combining this relation with \eqref{vahr} we prove the assertion \eqref{U} with $\alpha_1=\min\{\alpha_4;\alpha_5\}$.

To prove the additional assertion, let us consider the problem
\beq\label{tep9}
&&\hskip -1truecm u_1(t,x)=J_t^\beta (A(x)+q(t,x))u_1(t,x)+\phi_1(t,x)\, ,\quad x\in \Omega\, ,\; t\in (0,T),
\\ \label{tep10}
&&{\cal B}u_1(t,x)=0\, ,\;\; (t,x)\in
(0,T)\times {\partial\Omega}\,
\eeq
with $\phi_1=J_t^\beta (q_tu+\psi_t)$. Due to the proven properties of $u$ and the additional assumptions of $\psi,q$ we have
$q_tu+\psi_t\in L_p((0,T);L_p(\Omega))$, hence $\phi_1\in H_p^\beta([0,T];L_p(\Omega))$, $\phi_1(0,x)=0$.
Lemma 3 implies that the problem \eqref{tep9}, \eqref{tep10} has a unique solution $u_1\in U_p$ and $u_1(0,\cdot)=0$.
The next aim is to show that $u_t=u_1$. This will complete the proof.

Taking the operator $J_t^1$ from the relations
\eqref{tep9}, \eqref{tep10} and subtracing from \eqref{tep1}, \eqref{tep2} we see that $w=u-J_t^1u_1$ satisfies the following problem:
\beq\label{tep11}
&&\hskip -1truecm w(t,x)=J_t^\beta (A(x)+q(t,x))w(t,x)-J_t^{1+\beta}(q_t w)\, ,\quad x\in \Omega\, ,\; t\in (0,T),
\\ \label{tep12}
&&{\cal B}w(t,x)=0\, ,\;\; (t,x)\in
(0,T)\times {\partial\Omega}\, .
\eeq
Define $({\cal P}_t z)(\tau,\cdot)=\left\{\!\!\begin{array}{ll}
z(\tau,\cdot) &\mbox{if $\tau\le t$}\\ 0 &\mbox{if $\tau> t$}.\end{array}\right.$ Evidently, the relation
$({\mathcal S}J_t^\beta z)(\tau,\cdot)=\break ({\mathcal S}J_t^\beta {\cal P}_t z)(\tau,\cdot)$ is valid for $\tau\le t$
 and any $z\in L_p((0,T);L_p(\Omega))$, where ${\mathcal S}$ is the operator
defined in Lemma 3. Let us fix $t\in (0,T)$. Applying Lemma 3 to the problem \eqref{tep11}, \eqref{tep12} and using the embedding theorem we deduce the estimate
\beq\nonumber
&&\max_{0\le \tau\le t} \|w(\tau,\cdot)\|_{L_p(\Omega)}
=\max_{0\le \tau\le t}\|{\mathcal S}J_t^\beta{\cal P}_t J_t^{1}(q_t w)(\tau,\cdot)\|_{L_p(\Omega)}
\\ \nonumber
&&\quad
\le \tilde C_1\|{\mathcal S}J_t^\beta{\cal P}_t J_t^{1}(q_t w)\|_{U_p}\le \tilde C_1\|{\mathcal S}\|\, \|J_t^\beta{\cal P}_t J_t^{1}(q_t w)\|_{H_p^\beta([0,T];L_p(\Omega))}
\\ \nonumber
&&\quad
\le  \tilde C_2 \|{\cal P}_t J_t^{1}(q_t w)\|_{L_p((0,T);L_p(\Omega))}= \tilde C_2 \|J_t^{1}(q_t w)\|_{L_p((0,t);L_p(\Omega))}
\\ \nonumber
&&\quad
\le \tilde C_2\Bigl[\int_0^t \Bigl[\int_0^\tau \|(q_t w)(s,\cdot)\|_{L_p(\Omega)}ds\Bigl]^pd\tau\Bigl]^{1\over p}
\\ \label{contrest}
&&\quad
\le\tilde C_2 \|q_t\|_{L_{1}((0,T);L_\infty(\Omega))}\Bigl[\int_0^t \Bigl(\max_{0\le s\le \tau} \|w(s,\cdot)\|_{L_p(\Omega)}\Bigl)^p d\tau\Bigl]
^{1\over p},
\eeq
where $\tilde C_1$ and $\tilde C_2$ are some constants. Thus, we have proved the inequality\break
$\Bigl(\max\limits_{0\le \tau\le t} \|w(\tau,\cdot)\|_{L_p(\Omega)}\Bigl)^p
- \tilde C_3\int_0^t \Bigl(\max\limits_{0\le s\le \tau} \|w(s,\cdot)\|_{L_p(\Omega)}\Bigl)^p d\tau\le 0$ for
any $t\in (0,T)$ with some constant  $\tilde C_3$.
Making use of the the Gronwall's theorem we obtain\break $\max\limits_{0\le \tau\le t} \|w(\tau,\cdot)\|_{L_p(\Omega)}=0$ for any $t\in (0,T)$. This yields
$w=u-J_t^1u_1=0$ and in turn the desired equality $u_t=u_1$.
Theorem 1 is proved. \hfill $\Box$

\section{Positivity principle}\label{s:2a}
\setcounter{equation}{0}

In this section we prove a positivity principle for the solution of the equation \eqref{maineq} in a bit more general form. Namely, we consider the equation
\beq\label{maineqgen} D_t^{\{k\}} [u(t,x)-u_0(x)]\, =\, A(x)u(t,x) +
f(u(t,x),t,x)\, ,\quad t\in (0,T),\; x\in \Omega\, ,
 \eeq
where
$$
D_t^{\{k\}} v={d\over dt}k*v\, ,\quad (k*v)(t)=\int_0^t k(t-\tau)v(\tau)d\tau
$$
and the kernel $k$ has the following properties:
\beq\label{kassum}
&&k\in L_1(0,T)\cap C(0,T],\;\; k> 0,\, k - \mbox{decreasing} ,
\\ \label{kassum1}
&&k(t)\to\infty \quad \mbox{as}\quad t\to 0^+.
\eeq
This generalization doesn't make the proofs more complicated. Moreover, some parts of proofs even require
such a more general treatment.

\medskip
\noindent
{\bf Theorem 2}.
{\it Assume \eqref{kassum}, \eqref{kassum1},
$a_{ij},a_j\in C(\overline\Omega)$, $\omega\in (C({\partial\Omega}))^n$,
\beq
\label{as2}
&&\hskip -1truecm f\in C(\R\times [0,T]\times\overline\Omega)\, ,\;\;
\\
\label{as3}
&&\hskip -1truecm  \exists M\ge 0,\,\eta> 0\, :\,f(w,t,x)\ge -M|w|\quad\mbox{in $(-\eta,0)\times[0,T]\times\overline\Omega$}.
\eeq
Let $u\in C([0,T]\times \overline\Omega)$ with $u_{x_j},u_{x_i,x_j}\in C((0,T]\times \overline\Omega)$,
\beq\label{as3aa}
D_t^{\{k\}}(u-u_0)\in C((0,T]\times\overline\Omega)
\eeq
solves the problem \eqref{maineqgen}, \eqref{ini}, \eqref{bound}.
Moreover, let
\beq\label{as4}
\lim_{\epsilon\to 0^+}{1\over\epsilon}\int_0^\epsilon k(\tau)d\tau \sup_{0\le s\le \epsilon}|u(t-s,x)-u(t,x)|=0\quad \forall t\in (0,T],\, x\in\overline\Omega.
\eeq
Finally, we assume  $u_0\ge 0$ and $g\ge 0$.
Then the following assertions are valid:
\begin{description}
\item{\rm (i)} $u\ge 0$;
\item{\rm (ii)} if
 $u(t_0,x_0)=0$ at some point  $(t_0,x_0)\in (0,T]\times \Omega_N$, where $$\Omega_N=\left\{\begin{array}{ll}
\Omega &\mbox{in case {\rm I}}\\
\overline\Omega &\mbox{in case {\rm II}},\end{array}\right.$$ then $u(t,x_0)=0$ for any $t\in [0,t_0]$.
\end{description}
Recall that
the cases I and II were defined in \eqref{Bdef}.
}

\noindent
{\bf Remark 1}.
In case $k$ is the kernel of the fractional derivative $D_t^\beta$, i.e.  $k(t)={t^{-\beta}\over \Gamma(1-\beta)}$, $0<\beta<1$,
 the assumptions \eqref{as3aa}, \eqref{as4} are satisfied provided
$u\in C^{\beta'}([0,T];C(\overline\Omega))$ with some $\beta'>\beta$.
\\[1ex]
Before proving Theorem 2, we state and prove a lemma.
\\[1ex]
{\bf Lemma 5} (a minimum principle). {\it Let the assumptions of
Theorem {\rm 2} be satisfied, except for the conditions
\eqref{kassum1} and $u_0,g\ge 0$. Moreover, let instead of
\eqref{as3} the following stronger condition \beq\label{astr}
f(w,t,x)\ge 0\quad\mbox{in $(-\infty,0)\times
[0,T]\times\overline\Omega$} \eeq be satisfied. Let $(t_1,x_1)$
be the minimum
point of a solution of \eqref{maineqgen}, \eqref{ini} over
$[0,T]\times \overline\Omega$. If this
 is a stationary minimum with respect to $x$, i.e. $\nabla u(t_1,x_1)=0$, then
$u(t_1,x_1)\ge \min\limits_{x\in \overline\Omega}\{0;u_0(x)\}$.}
\\[1ex]
{\it Proof}.
Suppose that the assertion of the lemma doesn't hold.
Then
$$(t_1,x_1)\in(0,T]\times \overline\Omega \quad \mbox{and} \quad u(t_1,x_1)<0,\; u(t_1,x_1)<\min\limits_{x\in \overline\Omega}u_0(x).
$$
Since $x=x_1$ is the stationary minimum point of $u(t_1,x)$ over $\overline\Omega$ and the principal part of $A$ is elliptic, it holds
$
\sum\limits_{i,j=1}^na_{ij}(x_1)u_{x_ix_j}(t_1,x_1)\geq 0
$ (see \cite{Mir}).
Thus,
\begin{equation}\label{eq:47}
A(x_1)u(t_1,x_1)\ge 0.
\end{equation}
Further, in view of \eqref{astr} and $u(t_1,x_1)< 0$  we get
\beq\label{t1p3}
f(u(t_1,x_1),t_1,x_1)\ge 0.
\eeq

Next let us study the term $D_t^{\{k\}}[u(t,x)-u_0(x)]$ at $t=t_1$ and $x=x_1$.
We have
\beq \nonumber
&&D_t^{\{k\}}[u(t,x_1)-u_0(x_1)]\Bigl|_{t=t_1}
={\partial\over\partial t}k*[u(t,x_1)-u_0(x_1)]\Bigl|_{t=t_1}
\\ \nonumber
&&=\lim_{\epsilon\to 0+}{1\over\epsilon}\Bigl\{\int_0^{t_1} k(t_1-\tau)[u(\tau,x_1)-u_0(x_1)]d\tau
\\ \label{eqint1}
&&\qquad-\int_0^{t_1-\epsilon} k(t_1-\epsilon-\tau)[u(\tau,x_1)-u_0(x_1)]d\tau\Bigl\}.
\eeq
Observing that  $u(\tau,x_1)\ge u(t_1,x_1)$ holds for $0\le\tau\le t_1$ and   taking the relation \eqref{kassum} into account
we estimate the term between the brackets $\{\}$ in \eqref{eqint1}:
\beqst
&&\int_0^{t_1} k(t_1-\tau)[u(\tau,x_1)-u_0(x_1)]d\tau
\\
&&-\int_0^{t_1-\epsilon} k(t_1-\epsilon-\tau)[u(\tau,x_1)-u_0(x_1)]d\tau
\\
&&=\int_0^{t_1-\epsilon} [k(t_1-\tau)-k(t_1-\epsilon-\tau)][u(\tau,x_1)-u_0(x_1)]d\tau
\\
&&+\int_{t_1-\epsilon}^{t_1} k(t_1-\tau)[u(\tau,x_1)-u_0(x_1)]d\tau
\\
&&\le\int_0^{t_1-\epsilon} [k(t_1-\tau)-k(t_1-\epsilon-\tau)]d\tau\, [u(t_1,x_1)-u_0(x_1)]
\\
&&+\int_{t_1-\epsilon}^{t_1} k(t_1-\tau)[u(\tau,x_1)-u_0(x_1)]d\tau
\\
&&=\int_{t_1-\epsilon}^{t_1}k(\tau)d\tau\, [u(t_1,x_1)-u_0(x_1)]
+\int_{0}^{\epsilon} k(\tau)[u(t_1-\tau,x_1)-u(t_1,x_1)]d\tau.
\eeqst
Thus, from \eqref{eqint1} due to  \eqref{as4}, \eqref{kassum} and the inequalities $t_1>0$, $u(t_1,x_1)<u_0(x_1)$ we obtain
\beq \nonumber
&&D_t^{\{k\}}[u(t,x_1)-u_0(x_1)]\Bigl|_{t=t_1} \le \lim_{\epsilon\to 0+}
{1\over\epsilon}\Bigl[\int_{t_1-\epsilon}^{t_1}k(\tau)d\tau\, [u(t_1,x_1)-u_0(x_1)]
\\ \nonumber
&&\qquad+\int_{0}^{\epsilon} k(\tau)[u(t_1-\tau,x_1)-u(t_1,x_1)]d\tau\Bigl]
\\ \label{t1p2}
&&\qquad=k(t_1)\, [u(t_1,x_1)-u_0(x_1)]<0.
\eeq

The inequalities \eqref{t1p2}, \eqref{eq:47} and \eqref{t1p3} show that the left-hand side of the equation \eqref{maineqgen} is negative but the right-hand side
is nonnegative at $t=t_1$, $x=x_1$. This is a contradiction. Therefore, the assertion of the lemma is valid.
\hfill $\Box$
\\[1ex]
\noindent
{\it Proof of Theorem} 2. Firstly, let us prove (i) in case \eqref{as3} is replaced by the stronger condition \eqref{astr}.
Let again $(t_1,x_1)$ be the minimum point of $u$ over $[0,T]\times\overline\Omega$.
Suppose that (i) doesn't hold. Then $u(t_1,x_1)$ is negative.
In case I the point  $(t_1,x_1)$
is contained in the subset $(0,T]\times \Omega$. This implies that
$\nabla u(t_1,x_1)=0$ and in view of Lemma 5 we reach the contradiction. In case II the minimum point $(t_1,x_1)$ is contained in $(0,T]\times \overline\Omega$.
If $x_1\in\Omega$, we reach the contradiction in a similar manner.

Thus, it remains to consider the case II with $x_1\in {\partial\Omega}$.
In the minimum point at the boundary the inequality
$\omega(x_1)\cdot \nabla\, u(t_1,x_1)\le 0$ holds. This together
with  the assumption $\omega\cdot \nabla\, u|_{x\in {\partial\Omega}}=g\ge 0$ implies $\omega(x_1)\cdot \nabla\,
u(t_1,x_1)=0$. Moreover, $\tau\cdot \nabla\, u(t_1,x_1)=0$, where
$\tau$ is any tangential direction on ${\partial\Omega}$ at $x_1$, because
$x=x_1$ is the minimum point of the $x$-dependent function
$u(t_1,x)$ over the set ${\partial\Omega}$. Putting these relations together we see that
$\nabla u(t_1,x_1)=0$. Again, we reach the contradiction with Lemma 5. The assertion (i) is proved in case
\eqref{astr}.

\smallskip
Secondly, we prove (i) in the general case \eqref{as3}. Note that due to the continuity of $f$, the relation \eqref{as3} can immediately be extended to
$f(w,t,x)\ge -\Bigl(M+{D_q\over\eta}\Bigl)|w|$ in $(-q,0)\times[0,T]\times\overline\Omega$ where $D_q=\max\limits_{-q\le w\le 0\atop (t,x)\in [0,T]\times \overline\Omega}|f(w,t,x)|$
and $q>0$ is an arbitrary number. Let us define the new nonlinearity function:
$$
\widehat f(w,t,x)=\left\{\begin{array}{ll}
f(w,t,x)\; &\mbox{in case}\;\; w\ge -Q\\
f(-Q,t,x)\; &\mbox{in case}\;\; w< -Q,\end{array}\right.\quad \mbox{where}\;\; Q=\max_{t\in [0,T]\atop x\in \overline\Omega}|u(t,x)|.
$$
Then
\beq\label{widehatfting}
\widehat f(w,t,x)\ge -\widehat M|w| \quad \mbox{in $(-\infty,0)\times[0,T]\times\overline\Omega$},
\eeq
where $\widehat M=M+{D_Q\over\eta}$. According to the definition of $\widehat f$, we can rewrite the equation for $u$ \eqref{maineqgen}
in the following form:
\beq\label{maineqgenhat} D_t^{\{k\}} [u(t,x)-u_0(x)]\, =\, A(x)u(t,x) +
\widehat f(u(t,x),t,x),\; t\in (0,T),\, x\in \Omega.
 \eeq

Define $\tilde u(t,x)=e^{-\sigma t}u(t,x)$, where $\sigma>0$ is a number that we will specify later.
Then $\tilde u$ solves the following equation:
\beq\label{maineq2}
D_t^{\{\tilde k\}}[\tilde u(t,x)-u_0(x)] = A(x)\tilde u(t,x) + \tilde f(\tilde u(t,x),t,x) ,\, x\in \Omega,\, t\in (0,T),
\eeq
and satisfies the initial and boundary conditions
\beq\label{inikaks}
\tilde u(0,x)=u_0(x), \;x\in\Omega,\qquad {\cal B}\tilde u(t,x)=\tilde g(t,x),\;(t,x)\in (0,T)\times{\partial\Omega},
\eeq
where
\beqst
&&\hskip -5truemm \tilde k(t)=e^{-\sigma t}k(t)+\sigma\int_T^te^{-\sigma s} k(s) ds,\quad \tilde g(t)=e^{-\sigma t}g(t),
\\
&&\hskip -5truemm \tilde f(w,t,x)=e^{-\sigma t}\widehat f(e^{\sigma t}w,t,x)-\sigma w\int_0^T e^{-\sigma s}k(s)ds + \sigma u_0(x)\int_t^T e^{-\sigma s}k(s)ds.
\eeqst

We are going to show that the assumptions of the theorem are satisfied for the problem \eqref{maineq2},
\eqref{inikaks}.
The smoothness conditions $\tilde k\in L_1(0,T)\cap C(0,T]$, $\tilde f\in C(\R\times
[0,T]\times\overline \Omega)$, $\tilde u\in C([0,T]\times \overline\Omega)$, $\tilde u_{x_j},\tilde u_{x_i,x_j}\in C((0,T]\times \overline\Omega)$ and $D_t^{\{\tilde k\}}(\tilde u-u_0)\in C((0,T]\times\overline\Omega)$ directly follow from the
assumptions imposed on $k$, $f$ and $u$. Moreover, \eqref{as4} holds with $k$ and $u$  replaced by $\tilde k$ and $\tilde u$.
Observing the definition of $\tilde k$ and
the assumptions \eqref{kassum} we immediately see that  $\tilde k$ is decreasing.
Representing $\tilde k$ in the form $\tilde k (t)= e^{-\sigma T}k(t)+ \sigma\int_t^T e^{-\sigma s}(k(t)-k(s))ds$ and
using again \eqref{kassum} we prove the inequality $\tilde k> 0$.

Furthermore, using \eqref{widehatfting},  the relation $\sigma\int_0^T e^{-\sigma s}k(s) ds \to\infty$ as $\sigma\to\infty$, following from \eqref{kassum1}, as well as the inequalities $u_0\ge 0$, $k> 0$, we can show that there exists a sufficiently large
 $\sigma$ such that the condition \eqref{astr} is valid with $f$ replaced by $\tilde f$.

Now we see that we can apply the first part of the proof of (i) to the problem \eqref{maineq2}, \eqref{inikaks} to obtain
$\tilde u\ge 0$.\footnote{This is the point where we need results for more general kernel  $k(t)$ instead of
${t^{-\beta}\over\Gamma(1-\beta)}$.}
This
implies $u\ge 0$. The assertion (i) is completely proved.

\smallskip
Finally, we prove (ii). Suppose that this assertion does not hold. Then, due to the continuity of $u$, for some $(t_0,x_0)\in (0,T]\times \Omega_N$, such that
 $u(t_0,x_0)=0$,  it holds
\beq\label{t2p1}
\exists \delta>0,\, t_2,t_3\in (0,t_0),\, t_2<t_3\, :\,u(t,x_0)\ge \delta\;\;\mbox{for}\;\; t\in (t_2,t_3).
\eeq
We have
\beq\nonumber
&&D_t^{\{k\}}[u(t,x_0)-u_0(x)]\Bigl|_{t=t_0}=\lim_{\epsilon\to 0+}{1\over\epsilon}\Bigl\{\int_0^{t_0} k(t_0-\tau)[u(\tau,x_0)-u_0(x)]d\tau
\\ \label{ii1}
&&-\int_0^{t_0-\epsilon} k(t_0-\epsilon-\tau)[u(\tau,x_0)-u_0(x)]d\tau\Bigl\}.
\eeq
Let  $\epsilon<\min\{t_0-t_3;t_2\}$. Using \eqref{kassum}, \eqref{t2p1} and the relation $u\ge 0$ we estimate
the term between the brackets $\{\}$ in \eqref{ii1}:
\beq\nonumber
&&\int_0^{t_0} k(t_0-\tau)[u(\tau,x_0)-u(x_0)]d\tau
-\int_0^{t_0-\epsilon} k(t_0-\epsilon-\tau)[u(\tau,x_0)-u_0(x)]d\tau
\\ \nonumber
&&=\int_0^{t_0-\epsilon} [k(t_0-\tau)-k(t_0-\epsilon-\tau)]u(\tau,x_0)d\tau
\\ \nonumber
&&+\int_{t_0-\epsilon}^{t_0} k(t_0-\tau)u(\tau,x_0)d\tau-\int_{t_0-\epsilon}^{t_0} k(\tau) d\tau\, u_0(x)
\\ \nonumber
&&\le {\delta}\int_{t_2}^{t_3} [k(t_0-\tau)-k(t_0-\epsilon-\tau)]d\tau +
 \int_{t_0-\epsilon}^{t_0} k(t_0-\tau)u(\tau,x_0)d\tau
 \\ \label{ii2}
&&={\delta}\left[\int_{t_0-t_3}^{t_0-t_2}k(\tau)d\tau -
\int_{t_0-t_3-\epsilon}^{t_0-t_2-\epsilon}k(\tau)d\tau \right]+\int_0^\epsilon k(\tau)u(t_0-\tau,x_0)d\tau.
\eeq
In view of \eqref{as4} and the relation $u(t_0,x_0)=0$ it holds
 $\lim\limits_{\epsilon\to 0+}{1\over\epsilon}\int_0^\epsilon k(\tau)u(t_0-\tau,x_0)d\tau\break =0$.
Therefore, from \eqref{ii1} and \eqref{ii2} we get
\beqst
&&D_t^{\{k\}}[u(t,x_0)-u_0(x)]\Bigl|_{t=t_0}\le \lim_{\epsilon\to 0+}
{\delta\over\epsilon}\left[\int_{t_0-t_3}^{t_0-t_2}k(\tau)d\tau -
\int_{t_0-t_3-\epsilon}^{t_0-t_2-\epsilon}k(\tau)d\tau \right]
\\
&&=\delta {d\over ds}\int_{t_0-t_3+s}^{t_0-t_2+s}k(\tau)d\tau\Bigl|_{s=0}
=\delta[k(t_0-t_2)-k(t_0-t_3)]<0,
\eeqst
because $k$ is decreasing.

On the other hand,
 $(t_0,x_0)$ is a stationary local minimum point of $u(t_0,x)$, i.e. it holds $\nabla u(t_0,x_0)=0$. This follows from the assumption $u(t_0,x_0)=0$ and the inequality $u(t_0,x)\ge 0$ that holds
  in the neighborhood of $x_0$ as well as from the condition $\omega(x_0)\cdot \nabla u(t_0,x_0)=g(x_0)\ge 0$ in case II if $x_0\in{\partial\Omega}$.
Thus, $A(x_0)u(t_0,x_0)\ge 0$. Moreover, \eqref{as3} with \eqref{as2} implies $f(u(t_0,x_0),t_0,x_0)=u(0,t_0,x_0)\ge 0$. Consequently, the left-hand side of \eqref{maineqgen} is negative
but the right-hand side is nonnegative at $t=t_0$, $x=x_0$. This is a contradiction. The assertion (ii) is valid.
Theorem is completely proved. \hfill $\Box$
\\[2ex]
{\bf Corollary 1}. {\it Let the assumptions of Theorem {\rm 2} be satisfied. If $f(0,t_0,x_0)>0$ for some
$(t_0,x_0)\in (0,T]\times \Omega_N$, then $u(t_0,x_0)>0$.}
\\[1ex]
{\it Proof}. Due to Theorem 2, (i) we have $u\ge 0$. Assume that $f(0,t_0,x_0)>0$ for some
$(t_0,x_0)\in (0,T]\times \Omega_N$ and suppose contrary that $u(t_0,x_0)=0$.
Let us consider the equation \eqref{maineqgen} at  $t=t_0$, $x=x_0$. By Theorem 2 (ii) it holds
$u(t_0,x_0)=0$ for all $t\in [0,t_0]$, thus the left-hand side of this equation is zero. Since
$(t_0,x_0)$ is a stationary local minimum point of $u(t,\cdot)$, it holds $A(x_0)u(t_0,x_0)\ge 0$. Thus, in view of the assumption
$f(0,t_0,x_0)>0$, the right-hand side of the equation is positive. We reached a contradiction. The
assertion $u(t_0,x_0)>0$ is valid. \hfill $\Box$

\section{Uniqueness results  in the general case}\label{s:3}
\setcounter{equation}{0}

In this section we state and prove a uniqueness theorem for IP that is the most general in the sense of the Borel measure $\mu$.
But firstly we provide a technical lemma.

{\bf Lemma 6}. {\it
Let  $u\in C^\alpha ([0,T];L_s(\Omega))\cap C([0,T]\times\overline\Omega)$ with some $\alpha\in (0,1)$, $s\in [1,\infty]$ and
\beq\label{breveF}
f\in \breve{\cal  F}_{\alpha}:=C(\R;C^\alpha ([0,T];L_\infty(\Omega)))
\cap C^1(\R;L_\infty((0,T)\times\Omega)).
\eeq
Then
$f(u,\cdot,\cdot)\in C^\alpha([0,T];L_s(\Omega))$ and
\beq\label{le41}
\|f(u,\cdot,\cdot)\|_{C^\alpha([0,T];L_s(\Omega))}
\le \breve{\cal C}_{\hat m,\tilde m}\,\bigl(1+\|u\|_{C^\alpha([0,T];L_s(\Omega))}\bigl),
\eeq
where $\breve{\cal C}_{\hat m,\tilde m}$ is a constant depending on
$\hat m=\min\limits_{(t,x)\in [0,T]\times\overline\Omega}u(t,x)$ and $\tilde m=\break \max\limits_{(t,x)\in [0,T]\times\overline\Omega}u(t,x)$.
If, in addition,
$u\in C^\alpha([0,T];C^\gamma (\overline\Omega))$ with some $\gamma\in (0,1)$ and
\beqst
&&f\in {\cal F}_{\alpha,\gamma}:=C(\R;C^\alpha([0,T];C^\gamma(\overline\Omega)))
\\
&&\qquad\cap C^1\bigl(\R;C^\alpha([0,T];L_\infty(\Omega))\cap L_\infty((0,T);C^\gamma(\overline\Omega))\bigl)
\\
&&\qquad \cap C^2(\R;L_\infty((0,T)\times\Omega))
\eeqst
then $f(u,\cdot,\cdot)\in C^\alpha([0,T];C^\gamma (\overline\Omega))$ and
\beq\label{le42}
\|f(u,\cdot,\cdot)\|_{C^\alpha([0,T];C^\gamma (\overline\Omega))}\le {\cal C}_{\hat m,\tilde m}\,\bigl(1+\|u\|_{C^\alpha([0,T];C^\gamma (\overline\Omega))}\bigl)^2,
\eeq
where ${\cal C}_{\hat m,\tilde m}$ is a constant depending on $\hat m$ and $\tilde m$.}
\\[1ex]
{\it Proof}. We have $\|f(u,\cdot,\cdot)\|_{C^\alpha([0,T];L_s(\Omega))}=\sup\limits_{t\in (0,T)}\|f(u(t,\cdot),t,\cdot)\|_{L_s(\Omega)}+
\break \sup\limits_{t,\tau\in (0,T)}{1\over |t-\tau|^\alpha}\bigl\|f(u(s,\cdot),s,\cdot)|_{s=\tau}^{s=t}\,\bigl\|_{L_s(\Omega)}$.
Splitting the function under the second $\sup$ in this formula up in the following manner:
\beq\label{le43}
f(u(s,x),s,x)\bigl|_{s=\tau}^{s=t}\;=f(u(t,x),s,x)\bigl|_{s=\tau}^{s=t}+\int_{u(\tau,x)}^{u(t,x)}f_w(w,\tau,x)dw
\eeq
we deduce
\beqst
&&\|f(u,\cdot,\cdot)\|_{C^\alpha([0,T];L_s(\Omega))}\le \sup\limits_{w\in [\hat m,\tilde m]}\sup\limits_{t\in (0,T)}\|f(w,t,\cdot)\|_{L_\infty(\Omega)}({\rm meas}\,\Omega)^{1/s}
\\
&&\quad +\sup\limits_{w\in[\hat m,\tilde m]}\sup\limits_{t,\tau\in (0,T)}{1\over |t-\tau|^\alpha}\bigl\|f(w,s,\cdot)|_{s=\tau}^{s=t}\,\bigl\|_{L_\infty(\Omega)}({\rm meas}\,\Omega)^{1/s}
\\
&&\quad +\sup\limits_{w\in [\hat m,\tilde m]}\|f_w(w,\cdot,\cdot)\|_{L_\infty((0,T)\times\Omega)}
\sup\limits_{t,\tau\in (0,T)}{1\over |t-\tau|^\alpha}\bigl\|u(s,\cdot)|_{s=\tau}^{s=t}\,\bigl\|_{L_s(\Omega)}.
\eeqst
This implies $f(u,\cdot,\cdot)\in C^\alpha([0,T];L_s(\Omega))$ with \eqref{le41}.
Further, we have
\beqst
\begin{array}{ll}
&\|f(u(\cdot,\cdot),\cdot,\cdot)\|_{C^\alpha([0,T];C^\gamma(\overline\Omega))}=
\sup\limits_{t\in (0,T)}\Bigl[\,\sup\limits_{x\in\Omega}|f(u(t,x),t,x)|
\\
&+\sup\limits_{x,y\in\Omega}{\bigl|f(u(t,z),t,z)|_{z=y}^{z=x}\,\bigl|\over |x-y|^\gamma}\Bigl]
+\sup\limits_{t,\tau\in (0,T)}{1\over |t-\tau|^\alpha}\Bigl[\,\sup\limits_{x\in\Omega}\Bigl|f(u(s,x),s,x)\bigl|_{s=\tau}^{s=t}\,\Bigl|
\\
&+
\sup\limits_{x,y\in\Omega}{\bigl|f(u(s,z),s,z)|_{s=\tau}^{s=t}|_{z=y}^{z=x}\;\bigl|\over |x-y|^\gamma}\Bigl].
\end{array}
\eeqst
Splitting the terms in this formula up by means of the relations \eqref{le43} and the expressions
\beqst
&&f(u(t,z),t,z)\Bigl|_{z=y}^{z=x}=f(u(t,x),t,z)\Bigl|_{z=y}^{z=x}+\int_{u(t,y)}^{u(t,x)}f_w(w,t,y)dw,
\\
&&f(u(s,z),s,z)\Bigl|_{s=\tau}^{s=t}\Bigl|_{z=y}^{z=x}\;=f(u(t,x),s,z)\Bigl|_{s=\tau}^{s=t}\Bigl|_{z=y}^{z=x}\;
\\
&&+\int_{u(t,y)}^{u(t,x)} f_w(w,s,y)\Bigl|_{s=\tau}^{s=t}dw+\int_{u(\tau,x)}^{u(t,x)} f_w(w,\tau,z)\Bigl|_{z=y}^{z=x}dw
\\
&&+\int_0^{u(t,x)-u(\tau,x)}\int_{v+u(\tau,y)}^{v+u(\tau,x)} f_{ww}(w,\tau,y)dw\, dv
\\
&&+\int_{u(t,y)-u(\tau,y)}^{u(t,x)-u(\tau,x)}f_w(w+u(\tau,y),\tau,y)dw,
\eeqst
and
estimating we obtain $f(u,\cdot,\cdot)\in C^\alpha([0,T];C^\gamma (\overline\Omega))$ with \eqref{le42}.
\hfill $\Box$
\\[2ex]
{\bf Theorem 3}. {\it
Suppose that} IP {\it has two solutions $(z^j,u^j)$, $j=1,2$, such that
\beq\label{ipasu1}
z^j\in C^\alpha(\overline\Omega),\quad u^j\in C^{\alpha+\beta}([0,T];C^\alpha(\overline\Omega)),\; u^j_t\in C^{\bar\alpha}([0,T];L_\infty(\Omega))
\eeq
with some $\alpha\in (0,1)$ and $\bar\alpha\in (\max\{2\beta-1;0\},1)$.
Assume that $a_{ij}$, $a_j$, $\omega$ satisfy the conditions of Theorem {\rm 1} and
\beq
\label{ipas1}
\hskip -3truemm a,\, a_w, \, b_w\in {\cal F}_{\alpha+\beta,\alpha},\;
a_t,\, a_{ww},\, a_{wt}, \, b_{ww},\, b_{wt}\in \breve{\cal F}_{\bar\alpha}.
\eeq
 Moreover, let
\beq
\label{ipas4}
&&
\begin{array}{ll}
&\hskip -12truemm a_{ww}(w,t,x)z^2(x)+b_{ww}(w,t,x),\, a_{wt}(w,t,x)z^2(x)+b_{wt}(w,t,x)\ge 0
\\[.5ex]
&\quad\mbox{in}\;\;[\hat m,\tilde m]\times (0,T)\times\Omega,
\end{array}
\eeq
where $\hat m= \min\limits_{j=1,2}\min\limits_{(t,x)\in [0,T]\times\overline\Omega}u^j(t,x)$,
$\tilde m= \max\limits_{j=1,2}\max\limits_{(t,x)\in [0,T]\times\overline\Omega}u^j(t,x)$,
\beq \label{asumut}
&&u_t^1,u_t^2\ge 0 \quad\mbox{in case $a_{ww}\ne 0$ or $b_{ww}\ne 0$},
\\ \label{ipas2}
&&a(u_0,0,\cdot)=0,
\\ \label{ipas5}
&&a(u^1,\cdot,\cdot)\ge 0\, ,\;\;
\\ \label{ipas6}
&&\forall x\in \Omega \;\;\;\exists \varepsilon_x>0\, :\,
a(u^1(t,x),t,x)> 0\; \;\mbox{for}\;\;t\in (0,\varepsilon_x)
\eeq
and the relation
\beq
\label{ipas3a}
a_{w}(w,t,x)z^2(x)+b_{w}(w,t,x) \le \Theta \quad \mbox{in}\;\;[\hat m,\tilde m]\times (0,T)\times\Omega
\eeq
holds with
\beq \label{ipas52}
\Theta=\sup\{\theta\, :\, (D_t^\beta -\theta)a(u^1,\cdot,\cdot)\ge 0\}.
\eeq
Then $z^1=z^2$ and $u^1=u^2$. }
\\[1ex]
{\it Proof.}  Without restriction of generality assume that $\alpha\in (0,\min\{1-\beta;\bar\alpha+1-2\beta\})$ and
 $\bar\alpha\in (\max\{2\beta-1;0\},\beta)$.  Let us introduce the following notations:
\beq\label{ippqdef}
q(w,t,x)\, =\, a(w,t,x)z^2(x)+b(w,t,x)
\eeq
and $u=u^2-u^1$, $z=z^2-z^1$. Then $u$ satisfies the problem
\beq\label{ipp2}
&&\hskip -1truecm
\begin{split}
&D_t^\beta u(t,x)=A(x)u(t,x)+q(u^2(t,x),t,x)-q(u^1(t,x),t,x)
\\
&\qquad\qquad+a(u^1(t,x),t,x)z(x)\, ,\quad (t,x)\in (0,T)\times\Omega\,  ,
\end{split}
\\ \label{ipp3}
&&\hskip -1truecm u(0,x)=0,\, x\in\Omega\, ,\;\; {\cal B}u(t,x)=0\, ,\;\; (t,x)\in
(0,T)\times {\partial\Omega}.
\eeq
Further, define
$$
z^+={|z|+z\over 2}\, \;\;\; z^-={|z|-z\over 2}.
$$
Since $z\in C^\alpha(\overline\Omega)$, it hold $z^\pm \in C^\alpha(\overline\Omega)$. Now we consider the problems
\beq\label{ipp4}
&&\hskip -1truecm
\begin{split}
&D_t^\beta u^\pm(t,x)=A(x)u^\pm(t,x)+q^0(t,x)u^\pm(t,x)
\\
&\qquad\qquad+a(u^1(t,x),t,x)z^\pm(x)\, ,\quad (t,x)\in (0,T)\times\Omega\, ,\;
\end{split}
\\ \label{ipp5}
&&\hskip -1truecm u^\pm(0,x)=0,\, x\in\Omega\, ,\;\; {\cal B}u^\pm(t,x)=0\, ,\;\; (t,x)\in
(0,T)\times {\partial\Omega},
\eeq
where
\beq\label{ippq0def}
\begin{split}
&q^0(t,x)={q(u^2(t,x),t,x)-q(u^1(t,x),t,x)\over u^2(t,x)-u^1(t,x)}
\\
&\qquad =\int_0^1 q_w\bigl(u^1(t,x)(1-s)+u^2(t,x)s,t,x\bigl)ds.
\end{split}
\eeq
Due to the assumptions \eqref{ipasu1}, \eqref{ipas1} and Lemma 6 we have
\beq\label{q0prop}
&&q^0\in C^{\alpha+\beta}([0,T];C^\alpha(\overline\Omega))\, ,\quad  q^{0}_t\in C^{\bar\alpha}([0,T];L_\infty(\Omega)),
\\ \label{aprop}
&&\hskip-1.4truecm a(u^1,\cdot,\cdot)\in C^{\alpha+\beta}([0,T];\!C^\alpha(\overline\Omega)), {d\over dt}
a(u^1,\cdot,\cdot)\in
C^{\bar\alpha}([0,T];\!L_\infty(\Omega)).
\eeq
Moreover,
\beq\label{anull}
a(u^1(t,\cdot),t,\cdot)|_{t=0}=
a(u_0,0,\cdot)=0.
\eeq

Let $p$ be some number such that $\max\{{n\over 2}; {1\over\beta-\bar\alpha}\}<p<\infty$,
$p\not\in\{{1\over\beta}+{1\over 2};{2\over \beta}+1\}$.
By \eqref{q0prop}, \eqref{aprop}, \eqref{anull}  and  Theorem 1 we see that \eqref{ipp4}, \eqref{ipp5} have  solutions
 satisfying
\beq\label{upmspaces}
\begin{array}{ll}
&u^\pm\in C^{\alpha_1+\beta}([0,T];C^{\alpha_1}(\overline\Omega))\cap  C^{\alpha_1}([0,T];C^{2+\alpha_1}(\overline\Omega))
\;\;
\\[1ex]
&\mbox{with some
$\alpha_1\in (0,\alpha]$}\quad \mbox{and}\quad  u_t^\pm\in U_p.
\end{array}
\eeq
Moreover, taking the assumption
\eqref{ipas5}, the relations $z^\pm\ge 0$ and \eqref{ipp5} into account
and applying Theorem 2 with Remark 1  to the equation \eqref{ipp4},  we prove the inequality $u^\pm\ge 0$ and
the relation
\beq\label{ippfirst}
\begin{split}
&\mbox{if $u^{\rm s}(t_0,x_0)$ with some  ${\rm s}\in \{+;-\}$ equals zero in some point  }
\\
&\qquad \mbox{ $(t_0,x_0)\in (0,T]\times \Omega_N$ then $u^{\rm s}(t,x_0)=0$ for any $t\in [0,t_0]$.}
\end{split}
\eeq
Subtracting the problems for $u^-$ and $u$ from the problem for $u^+$ and observing that  $z=z^+ -z^-$, we see that the function
$\overline u=u^+-u^--u$ solves the homogeneous equation $D_t^\beta \overline u=(A+q^0)\overline u$ and satisfies the homogeneous
initial and boundary conditions. Therefore, due to the uniqueness assertion of Theorem 1 we get $\overline u=0$, hence
\beq\label{breakingpoint}
u=u^+-u^-\, .
\eeq

Next let us consider the problems
\beq\label{ipp6}
&&\hskip -1truecm
\begin{split}
&D_t^\beta v^\pm(t,x)=A(x)v^\pm(t,x)+q^0(t,x)v^\pm(t,x)
\\
&\qquad\qquad+\varphi^\pm(t,x)\, ,\quad (t,x)\in (0,T)\times\Omega\, ,
\end{split}
\\ \label{ipp7}
&&\hskip -1truecm v^\pm(0,x)=0,\, x\in\Omega\, ,\;\; {\cal B}v^\pm(t,x)=0\, ,\;\; (t,x)\in
(0,T)\times {\partial\Omega},
\eeq
where
$$
\varphi^\pm = D_t^\beta (q^0 u^\pm)-q^0 D_t^\beta u^\pm + (D_t^\beta -\Theta)a(u^1(t,x),t,x)z^\pm.
$$
Analyze the formula of $\varphi^\pm$ termwise.
In view of \eqref{q0prop} and $u^\pm\in C^{\alpha_1+\beta}([0,T];C^{\alpha_1}(\overline\Omega))$
it hold $q^0 u^\pm\in C^{\alpha_1+\beta}([0,T];C^{\alpha_1}(\overline\Omega))$. Moreover, $(q^0 u^\pm)(0,\cdot)=0$. Therefore, due to
Lemma 2 we have $D_t^\beta (q^0 u^\pm)\in C^{\alpha_1}([0,T];C^{\alpha_1}(\overline\Omega))$ and $D_t^\beta (q^0 u^\pm)|_{t=0}=0$.
Similarly, we deduce $q^0 D_t^\beta u^\pm
\in C^{\alpha_1}([0,T];C^{\alpha_1}(\overline\Omega))$, $q^0 D_t^\beta u^\pm|_{t=0}=0$. In virtue of \eqref{aprop}, \eqref{anull} and Lemma 2 we get $(D_t^\beta -\Theta)a(u^1,\cdot,\cdot)z^\pm \in C^{\alpha}([0,T];C^{\alpha}(\overline\Omega))$ and
$(D_t^\beta -\Theta)a(u^1(0,\cdot),0,\cdot)z^\pm=0$.
Thus, we see that $\varphi^\pm\in C^{\alpha_1}([0,T];C^{\alpha_1}(\overline\Omega))$ and
$\varphi^\pm (0,\cdot)=0$.

 Further, due to the range of $p$ and an embedding theorem it holds $U_p\subset C^{\bar\alpha}([0,T];\break L_p(\Omega))$.
 Thus, taking into account
  \eqref{q0prop} and \eqref{upmspaces}  we have
$(q^0u^\pm)_t\in\break C^{\bar\alpha}([0,T];L_p(\Omega))$. Thus, taking into account the ranges of $\alpha_1$, $\alpha$, $\bar\alpha$
and applying Lemma 2
we obtain $D_t^\beta (q^0 u^\pm)=
J_t^{1-\beta}[(q^0 u^\pm)_t]\in C^{\bar\alpha+1-\beta}([0,T];L_p(\Omega))\subset C^{\alpha_1+\beta}([0,T];L_p(\Omega))$.
Similarly we prove $q^0 D_t^\beta u^\pm\in C^{\alpha_1+\beta}([0,T];L_p(\Omega))$. Finally, by the relation
$(D_t^\beta-\Theta)a(u^1,\cdot,\cdot)z^\pm=(J_t^{1-\beta}{d\over dt}-\Theta)a(u^1,\cdot,\cdot)z^\pm$,
\eqref{aprop} and Lemma 2 it holds $(D_t^\beta-\Theta)a(u^1(\cdot,\cdot),\cdot,\cdot)z^\pm\in
C^{\bar\alpha+1-\beta}([0,T];L_\infty(\Omega))\subset C^{\alpha_1+\beta}([0,T];L_p(\Omega))$.
Summing up, $\varphi^\pm \in C^{\alpha_1+\beta}([0,T];L_p(\Omega))$.

We see that the assumptions of
 Theorem 1 are satisfied for the problems \eqref{ipp6}, \eqref{ipp7}. This implies that \eqref{ipp6}, \eqref{ipp7} have solutions
 $v^\pm\in C^{\alpha_2+\beta}([0,T];C^{\alpha_2}(\overline\Omega))\cap  C^{\alpha_2}([0,T];C^{2+\alpha_2}(\overline\Omega))$ with some
$\alpha_2\in (0,\alpha_1]$.

Taking into account  the smoothness properties of $q^0$ and $u^\pm$, the relation $u^\pm (0,\cdot)\break =0$ and performing some elementary
computations we reach the following expression:
$$
\bigl(D_t^\beta (q^0u^\pm)-q^0 D_t^\beta u^\pm\bigl)(t,x)=\int_0^t {\beta(t-\tau)^{-\beta-1}\over \Gamma(1-\beta)}\bigl[q^0(t,x)-
q^0(\tau,x)\bigl]u^\pm (\tau,x)d\tau.
$$
By the assumptions \eqref{ipas4}, \eqref{asumut} and the formulas \eqref{ippqdef}, \eqref{ippq0def} we see that the function
$q^0$ is nondecreasing in $t$. This with the proven inequalities $u^\pm\ge 0$ implies that $D_t^\beta (q^0u^\pm)-q^0 D_t^\beta u^\pm\ge 0$.
Moreover, by \eqref{ipas52} and $z^\pm\ge 0$ it hold $(D_t^\beta -\Theta)a(u^1(t,x),t,x)z^\pm\ge 0$. Consequently, we have
$\varphi^\pm\ge 0$.
Applying Theorem 2 with Remark 1  to the equation \eqref{ipp6},  we obtain the inequality $v^\pm\ge 0$ and
the relation
\beq\label{ippsec}
\begin{split}
&\mbox{if $v^{\rm s}(t_0,x_0)$ with some  ${\rm s}\in \{+;-\}$ equals zero in some point  }
\\
&\qquad \mbox{ $(t_0,x_0)\in (0,T]\times \Omega_N$ then $v^{\rm s}(t,x_0)=0$ for any $t\in [0,t_0]$.}
\end{split}
\eeq

Next we are going to establish  relations between $v^\pm$ and $u^\pm$. To this end, we deduce and analyze problems for the functions
$Q^\pm=u^\pm-J_t^\beta \bigl(v^\pm+\Theta u^\pm\bigl)$. Adding \eqref{ipp4} multiplied by $\Theta$ to \eqref{ipp6}, taking the operator
$J_t^\beta$ and subtracting from \eqref{ipp4} after some transformations we arrive at the following equation for $Q^\pm$:
\beq\label{ippQ}
\begin{split}
&D_t^\beta Q^\pm(t,x)=A(x)Q^\pm(t,x)+q^0(t,x) Q^\pm(t,x)+\zeta^\pm (t,x),
\\
&\; (t,x)\in (0,T)\times\Omega,
\end{split}
\eeq
with the initial and boundary conditions
\beq\label{ippQb}
Q^\pm(0,x)=0,\, x\in\Omega\, ,\;\; {\cal B}Q^\pm(t,x)=0\, ,\;\; (t,x)\in
(0,T)\times {\partial\Omega},
\eeq
where
\beqst
&&\zeta^\pm (t,x)=\bigl[J_t^\beta (q^0 D_t^\beta Q^\pm)-q^0 Q^\pm\bigl](t,x)
\\
&&\quad=\int_0^t {(t-\tau)^{\beta-1}\over \Gamma(\beta)}[q^0(\tau,x)-q^0(t,x)](D_t^\beta Q^\pm )(\tau,x)d\tau
\\
&&\quad
=\int_0^t \Bigl\{{(\beta-1)(t-\tau)^{\beta-2}\over \Gamma(\beta)}[q^0(\tau,x)-q^0(t,x)]
\\
&&\qquad\quad -{(t-\tau)^{\beta-1}\over \Gamma(\beta)}q^0_\tau(\tau,x)\Bigl\}(J_t^{1-\beta} Q^\pm) (\tau,x)d\tau.
\eeqst
Using \eqref{q0prop} and the relation $J_t^\beta J_t^{1-\beta}=J_t^1$ we deduce
\beq\nonumber
&&\|\zeta^\pm(t,\cdot) \|_{L_p(\Omega)}\le \hat C_1\int_0^t {(t-\tau)^{\beta-1}\over \Gamma(\beta)}
\int_0^\tau {(\tau-s)^{-\beta}\over \Gamma(1-\beta)}\|Q^\pm (s,\cdot\|_{L_p(\Omega)}ds d\tau
\\ \label{ippQc}
&&\qquad =\hat C_1 \int_0^t \|Q^\pm (\tau,\cdot\|_{L_p(\Omega)}d\tau
\eeq
with some constant $\hat C_1$. Estimating the solution of
\eqref{ippQ}, \eqref{ippQb} by means of the technique that was used in derivation of the estimate \eqref{contrest}
we obtain\break $\max\limits_{0\le \tau\le t}\|Q^\pm(\tau,\cdot) \|_{L_p(\Omega)}\le \hat C_2 \|\zeta^\pm\|_{L_p((0,t);L_p(\Omega))}$
for any $t\in (0,T)$ with some constant $\hat C_2$. This with \eqref{ippQc} implies
\beqst
&&\Bigl[\max_{0\le \tau\le t}\|Q^\pm(\tau,\cdot) \|_{L_p(\Omega)}\Bigl]^p\le \hat C_1 \hat C_2 T^p
\int_0^t \Bigl[\max_{0\le \tau\le s}\|Q^\pm(\tau,\cdot) \|_{L_p(\Omega)}\Bigl]^p ds
\eeqst
for any $t\in (0,T)$. Applying Gronwall's theorem we reach the equality\break
$\max\limits_{0\le \tau\le t}\|Q^\pm(\tau,\cdot) \|_{L_p(\Omega)}=0$ for any $t\in (0,T)$. Therefore, $Q=u^\pm-J_t^\beta \bigl(v^\pm+\Theta u^\pm\bigl)=0$.
This yields the following relations between $v^\pm$ and $u^\pm$:
\beq\label{ipp100}
u^\pm-\Theta J_t^\beta u^\pm=J_t^\beta v^\pm \quad \Leftrightarrow\quad  v^\pm=D_t^\beta u^\pm-\Theta u^\pm.
\eeq

Let us return to the equation \eqref{ipp4}. We subtract the term $\Theta u^\pm$ from both sides and integrate. Taking into account
the right relation in \eqref{ipp100} we obtain
\beq\label{ipp101}
&&\hskip -1truecm
\begin{split}
&\int_0^T v^\pm (t,x)d\mu=A(x)\int_0^T u^\pm (t,x)d\mu+\int_0^T\bigl[q^0(t,x)-\Theta\bigl]u^\pm(t,x)d\mu
\\
&\qquad\qquad+\int_0^Ta(u^1(t,x),t,x)d\mu\, z^\pm(x)
\end{split}
\eeq
for any $x\in \Omega$. By continuity, this relation can be extended to $\overline \Omega$. Due to the relations $u^2-u^1=u^+-u^-$ and
the equality $\int_0^T u^1(t,x)d\mu=\int_0^T u^2(t,x)d\mu$ for $x\in \Omega$, following from \eqref{add}, we have
\beq\label{ipp102}
\int_0^T u^+(t,x)d\mu=\int_0^T u^-(t,x)d\mu
\eeq
in $\Omega$. By continuity, \eqref{ipp102} holds in $x\in\overline\Omega$. Define
\beq\label{ipp103}
x^*={\rm arg}\max_{x\in\overline\Omega}\int_0^T u^\pm(t,x)d\mu.
\eeq
Without restriction of generality we may assume that $x^*\in\Omega_N$. Indeed, if\break  $x^*\in\overline\Omega\setminus\Omega_N={\partial\Omega}$ in case I, the vanishing boundary condition $u^\pm|_{\partial\Omega} =0$
implies $\int_0^T u^\pm(t,x^*)d\mu=0$ and since $u^\pm\ge 0$ from \eqref{ipp103} we get $\int_0^T u^\pm(t,x)d\mu\equiv 0$, which means that we can redefine $x^*$
as an arbitrary point in $\Omega$.

According to the definitions of $z^\pm$ it holds either $z^+(x^*)=0$ or $z^-(x^*)=0$. Suppose that  $z^+(x^*)=0$. Then from \eqref{ipp101}
we get
\beq\label{ipp104}
&&\hskip -1.5truecm
\int_0^T\!\! v^+ (t,x^*)d\mu=A(x^*)\!\int_0^T\! \!u^+ (t,x^*)d\mu+\!\int_0^T\!\bigl[q^0(t,x^*)-\Theta\bigl]u^+(t,x^*)d\mu.
\eeq
Let us analyze the signs of the terms in  \eqref{ipp104}. By virtue of the proved relation $u^+\ge 0$ and the
equality ${\cal B}\int_0^T u^+(\cdot,x)=0$, $x\in{\partial\Omega}$, following from \eqref{ipp5}, $x^*$ is a stationary maximum point of
the function $\int_0^T u^\pm(t,x)d\mu$. This implies that
$A(x^*)\!\int_0^T\! \!u^+ (t,x^*)d\mu\le 0$. Moreover, the assumption \eqref{ipas3a} with \eqref{ippqdef}, \eqref{ippq0def}
and $u^+\ge 0$
implies $\int_0^T\!\bigl[q^0(t,x^*)-\Theta\bigl]u^+(t,x^*)d\mu\le 0$. Consequently, the right-hand side of \eqref{ipp104} is non-positive. Thus
 $\int_0^T\!\! v^+ (t,x^*)d\mu\le 0$. This with the proven relation $v^+\ge 0$
 implies $\int_0^T\!\! v^+ (t,x^*)d\mu= 0$ and
$v^+(t,x^*)=0$ for any $t\in {\rm supp}(\mu)$. In particular, $v^+(t^*,x^*)=0$, where
$$t^*=\sup {\rm supp}(\mu).$$
It holds
$t^*>0$ because of the assumption ${\rm supp}(\mu)\cap (0,T]\ne \emptyset$. Observing the relation \eqref{ippsec} we deduce
$v^+(t,x^*)=0$ for any $t\in [0,t^*]$.

In view of the left equality in \eqref{ipp100} we arrive at the following homogeneous Volterra equation of the 2. kind:
$(u-\Theta J_t^\beta) u^+(t,x^*)=0,\;\; t\in [0,t^*]$.
It has only the trivial solution
\beq\label{lopualgus}
u^+(t,x^*)=0\, ,\quad t\in [0,t^*].
\eeq
 Recall that this result is based on the supposition that
$z^+(x^*)=0$. Similarly we reach the equality $u^-(t,x^*)=0$, $t\in [0,t^*]$ if we suppose $z^-(x^*)=0$. Consequently, (see also
\eqref{ipp102}) we obtain the  relation $\int_0^T u^\pm (t,x^*)d\mu=0$.

Recall that $x^*$ is the maximum point of the non-negative function $\int_0^T u^\pm(t,x)d\mu$.
Therefore, it holds $\int_0^T u^\pm(\cdot,x)d\mu=0$ for any $x\in\Omega$. Since $u^\pm\ge 0$, we get  $u^\pm(t,x)=0$ for any
$t\in {\rm supp}(\mu)$ and $x\in \Omega$, in particular $u^\pm(t^*,x)=0$ for any  $x\in \Omega$. Applying
\eqref{ippfirst} we deduce
$$
u^\pm(t,x)=0 \quad \mbox{for any  $t\in [0,t^*]$ and $x\in \Omega$}.
$$
Due to this relation, the equation \eqref{ipp4} reduces to the form
$$
a(u^1(t,x),t,x)z^\pm (x)=0
$$
in the cylinder $(t,x)\in (0,t^*)\times \Omega$. In view of the assumption \eqref{ipas6} we obtain $z^\pm(x)=0$ in $\Omega$. Observing
that $z^2-z^1=z=z^+-z^-$ we reach the assertion $z^1=z^2$. Finally, due to $z^\pm=0$, the linear problem \eqref{ipp4}, \eqref{ipp5} is homogeneous.
By virtue of the uniqueness statement of Theorem 1 we get $u^\pm=0$. This in view of $u^2-u^1=u=u^+-u^-$  yields $u^1=u^2$. Theorem is
completely proved.
\hfill $\Box$

\medskip
From the proved theorem we  infer the following uniqueness result for a linear inverse problem.
\\[1ex]
{\bf Corollary 2}. {\it Let $a(w,t,x)=a(t,x)$ and $b(w,t,x)=b^1(t,x)w+b^2(t,x)$.
Suppose that} IP {\it has two solutions $(z^j,u^j)$, $j=1,2$, satisfying \eqref{ipasu1}
with some $\alpha\in (0,1)$ and $\bar\alpha\in (\max\{2\beta-1;0\},1)$. Assume that $a_{ij}$, $a_j$, $\omega$ fulfill the conditions of Theorem {\rm 1},
$a,b^1\in C^{\alpha+\beta}([0,T];C^\alpha(\Omega))$, $a_t, b_t^1\in C^{\bar\alpha}([0,T],L_\infty(\Omega))$,
$a,b^1_t\ge 0$, $a(0,\cdot)=0$,
\beq\label{aposlinear}
\forall x\in \Omega \;\;\;\exists \varepsilon_x>0\, :\,
a(t,x)> 0\; \;\mbox{for}\;\;t\in (0,\varepsilon_x)
\eeq
and
\beq\label{bposlinear}
b^1\le \Theta=\sup\{\theta\, :\, (D_t^\beta -\theta)a\ge 0\}.
\eeq
Then $z^1=z^2$ and $u^1=u^2$.
}

\section{Uniqueness in a particular case}\label{s:4}
\setcounter{equation}{0}

In case the measure $\mu$ has a special form, the uniqueness can be proved under lower regularity assumptions and
the cone condition \eqref{ipas4} and the restriction \eqref{asumut} dropped.
Let us consider the following particular case:
\beq\label{muspecial}
\mbox{$d\mu= \varkappa(t)dt$, where $dt$ is the Lebesgue measure and $\varkappa\ge 0$, $\varkappa\ne 0$.}
\eeq

\noindent
{\bf Theorem 4}. {\it Let \eqref{muspecial} hold, where $\varkappa\in W_s^1(0,T)$ with some $s>{1\over 1-\beta}$.
Suppose that} {IP} {\it has two solutions $(z^j,u^j)$, $j=1,2$, such that
\beq\label{ipasu1spec}
z^j\in C^\alpha(\overline\Omega),\quad u^j\in C^{\alpha+\beta}([0,T];L_\infty(\Omega))\cap  C^{\alpha}([0,T];C^\alpha(\overline\Omega))
\eeq
with some $\alpha\in (0,1)$.
Assume that $a_{ij}$, $a_j$, $\omega$ satisfy the conditions of Theorem {\rm 1} and
\beq
\label{ipas1spec}
\hskip -3truemm a,  a_w, b_w\in \breve{\cal F}_{\alpha+\beta}\cap {\cal F}_{\alpha,\alpha}.
\eeq
Moreover, let the relations \eqref{ipas2}, \eqref{ipas5}, \eqref{ipas6} hold and
the inequality  \eqref{ipas3a} be valid with some $\Theta<\hat\theta=\sup\{\theta\, :\, \kappa-\theta\varkappa\ge 0\}$,
where
\beq\label{kappaspec}
\kappa(t)={1\over \Gamma(1-\beta)}\Bigl[(T-t)^{-\beta}\varkappa(T)-\int_t^T(\tau-t)^{-\beta}\varkappa'(\tau)d\tau\Bigl].
\eeq
Then $z^1=z^2$ and $u^1=u^2$.
}
\\[1ex]
{\it Proof.} The beginning of the proof coincides with the section of the proof of Theorem 3 that starts with the
formula \eqref{ippqdef} and ends with  the relation \eqref{ippq0def}.
From \eqref{ippq0def} on we continue in a different manner. Due
to assumptions of the theorem and Lemma 6, it hold
\beq\label{q0propsepc}
&&q^0, a(u^1,\cdot,\cdot)\in C^{\alpha+\beta}([0,T];L_\infty(\Omega))\cap C^{\alpha}([0,T];C^\alpha(\overline\Omega)).
\eeq
Taking also the relation  $a(u^1(t,\cdot),t,\cdot)|_{t=0}=
a(u_0,0,\cdot)=0$ into account and using Theorem 1 we come to the conclusion that \eqref{ipp4}, \eqref{ipp5} have unique solutions
 $u^\pm\in C^{\alpha_1+\beta}([0,T];C^{\alpha_1}(\overline\Omega))\cap  C^{\alpha_1}([0,T];C^{2+\alpha_1}(\overline\Omega))$ with some
$\alpha_1\in (0,1-\beta)$ and $D_t^\beta u^\pm\in C^{\alpha_1}([0,T];C^{\alpha_1}(\overline\Omega))$.
Like in the proof of Theorem 3, we deduce the relations $u^\pm\ge 0$, \eqref{ippfirst} and the equality \eqref{breakingpoint}.

Let us subtract the term $\Theta u^\pm$ from the left and right-hand side of \eqref{ipp4} and integrate \eqref{ipp4} in the Lebesgue sense with the weight $\varkappa$.
Integrating by parts at the
left-hand side and performing some simple computations we obtain
\beq\label{ipp101spec}
&&\hskip -.7truecm
\begin{split}
&\int_0^T [\kappa(t)-\Theta\varkappa(t)]u^\pm (t,x)dt
=A(x)\int_0^T u^\pm (t,x)\varkappa(t) dt
\\
&\qquad+\int_0^T[q^0(t,x)-\Theta]u^\pm(t,x)\varkappa(t) dt+\int_0^Ta(u^1(t,x),t,x)\varkappa(t) dt\, z^\pm(x)
\end{split}
\eeq
for any $x\in \overline\Omega$.  By the relations $u^2-u^1=u^+-u^-$,
the equality $\int_0^T u^1(t,x)\varkappa(t)dt=\int_0^T u^2(t,x)\varkappa(t)dt$ for $x\in \Omega$ and continuity we have
\beq\label{ipp102spec}
\int_0^T u^+(t,x)\varkappa(t)dt=\int_0^T u^-(t,x)\varkappa(t)dt
\eeq
in $\overline\Omega$. Define
\beq\label{ipp103u}
x^*={\rm arg}\max_{x\in\overline\Omega}\int_0^T u^\pm(t,x)\varkappa(t)dt.
\eeq
Without restriction of generality we may assume that $x^*\in\Omega_N$.

It holds either $z^+(x^*)=0$ or $z^-(x^*)=0$. Let $z^+(x^*)=0$. Then from \eqref{ipp101spec}
we get
\beq\label{ipp104spec}
\begin{split}
&\int_0^T [\kappa(t)-\Theta\varkappa(t)]u^+ (t,x^*)dt
=A(x)\int_0^T u^+ (t,x^*)\varkappa(t) dt
\\
&\qquad+\int_0^T[q^0(t,x^*)-\Theta]u^+(t,x^*)\varkappa(t) dt.
\end{split}
\eeq
Arguing similarly as in the case of the equation \eqref{ipp104}, we show that
 the right-hand side of \eqref{ipp104spec} is non-positive.
Thus, we obtain
\beq\label{speca}
\int_0^T [\kappa(t)-\Theta\varkappa(t)]u^+ (t,x^*)dt\le 0.
\eeq
According to the definitions of $\Theta$ and $\hat\theta$ we have $\kappa-\Theta\varkappa=\kappa-\hat\theta\varkappa+(\hat\theta-\Theta)\varkappa\ge (\hat\theta-\Theta)\varkappa$, where
$\hat\theta-\Theta>0$. Thus, $\kappa-\Theta\varkappa\ge 0$ and due to the
relation ${\rm supp}(\varkappa)\cap (0,T]\ne\emptyset$, following from the assumptions
imposed on $\varkappa$, there exists $t_1\in (0,T]$ such that $\kappa(t_1)-\Theta\varkappa(t_1)>0$. Therefore,
observing also the property  $u^+\ge0$, from \eqref{speca} we get $[\kappa-\Theta\varkappa]u^+(\cdot,x^*)=0$ and $u^+(t_1,x^*)=0$.
Now the relation \eqref{ippfirst} implies $u^+(t,x^*)=0$, $t\in [0,t_1]$.

The rest of the proof repeats the part of the proof of Theorem 3 that follows the
formula \eqref{lopualgus} (with $t_1$ instead of $t^*$).
\hfill $\Box$

\medskip\noindent
{\bf Corollary 3}. {\it Let \eqref{muspecial} hold, where $\varkappa\in W_s^1(0,T)$ with some $s>{1\over 1-\beta}$
and $a(w,t,x)=a(t,x)$, $b(w,t,x)=b^1(t,x)w+b^2(t,x)$.
Suppose that} IP {\it has two solutions $(z^j,u^j)$, $j=1,2$, satisfying \eqref{ipasu1spec}
with some $\alpha\in (0,1)$.
Assume that $a_{ij}$, $a_j$, $\omega$ satisfy the conditions of Theorem {\rm 1}, $a,b^1\in C^\alpha([0,T];C^\alpha(\overline\Omega))$,
$a\ge 0$, $a(0,\cdot)=0$, \eqref{aposlinear} is valid and
$b^1\le \Theta<\hat\theta$  with $\hat\theta$ defined in Theorem {\rm 4}. Then $z^1=z^2$ and $u^1=u^2$.
}

\medskip\noindent
{\bf Remark 2}. In case $\varkappa=1$ we have $\kappa(t)={(T-t)^{-\beta}\over \Gamma(1-\beta)}$
 and $\hat\theta={T^{-\beta}\over \Gamma(1-\beta)}$.

\section{Additional remarks}\label{s:5}
\setcounter{equation}{0}

The aim of this section is to interpret
  the conditions \eqref{asumut},  \eqref{ipas5}, \eqref{ipas6}  in a suitable way and
  estimate the quantity $\Theta$ occurring in the conditions \eqref{ipas3a}, \eqref{bposlinear} from below.


\smallskip
\noindent
{\bf Theorem 5}. {\it Let $u$ solve the problem \eqref{maineq}, \eqref{ini}, \eqref{bound} with $f(w,t,x)=a(w,t,x)z(x)\break +b(w,t,x)$. Assume
$u\in C([0,T]\times\overline\Omega)$,
$u,u_{x_j},u_{x_ix_j} \in C^{\beta''}([0,T];L_1(\Omega))$ with some $\beta''>1-\beta$, $u_t\in
C^{\beta'}([0,T];C(\overline\Omega))$ with some $\beta'>\beta$, $u_{tx_j},u_{tx_ix_j}\in C((0,T]\times\overline\Omega)$,
$a_{ij},a_j\in C(\overline\Omega)$, $\omega\in (C({\partial\Omega}))^n$, $a,b\in C(\R;C^{\beta''}([0,T];L_\infty(\Omega)))$,
$a_w,b_w,a_t,b_t\in C(\R\times [0,T]\times\overline\Omega)$,  $a(u_0,0,\cdot)=0$, $Au_0+b(u_0,0,\cdot)=0$, $g_t\ge 0$ and
\beq\label{utpossuff}
a_{t}(u,\cdot,\cdot)z+b_{t}(u,\cdot,\cdot)\ge 0.
\eeq
Then $u_t\ge 0$. If in addition
\beq\label{utpossuff2}
a(u,\cdot,\cdot),\, a_{w}(u,\cdot,\cdot),\,  a_{t}(u,\cdot,\cdot)\ge 0,
\eeq
then
\beq\label{utposs1}
\Theta=\sup\{\theta\, :\, (D_t^\beta -\theta)a(u,\cdot,\cdot)\ge 0\}\ge {T^{-\beta}\over \Gamma(1-\beta)}.
\eeq
}
\noindent
{\it Proof.}
Firstly, let us prove $u_t\ge 0$.  The equation for $u$ can be rewritten in the form $\int_0^t{(t-\tau)^{-\beta}\over \Gamma(1-\beta)} u_\tau(\tau,x)d\tau=A(x)u(t,x)+
a(u(t,x),t,x)z(x)+b(u(t,x),t,x)$. The left hand side of this equation equals $J^{1-\beta}_tu_t$. Therefore, applying the operator $D_t^{1-\beta}$ to this equation we get
$u_t=D_t^{1-\beta}[Au+a(u,\cdot,\cdot)z+b(u,\cdot,\cdot)]$. By virtue of the assumptions of the theorem and Lemma 6 we have
$Au+a(u,\cdot,\cdot)z+b(u,\cdot,\cdot)\in C^{\beta''}([0,T];L_1(\Omega))$
and $[Au+a(u,\cdot,\cdot)z+b(u,\cdot,\cdot)]\Bigl|_{t=0}=0$. Applying Lemma 2 we get $D_t^{1-\beta}[Au+a(u,\cdot,\cdot)z+b(u,\cdot,\cdot)]\Bigl|_{t=0}=0$. Thus,
$u_t(0,\cdot)=0$. Now from the equation of $u$ we deduce the following equation for $u_t$:
$D_t^\beta [u_t-u_t(0,\cdot)]=Au_t+[a_w(u,\cdot,\cdot)z+b_w(u,\cdot,\cdot)]u_t+a_t(u,\cdot,\cdot)z+b_t(u,\cdot,\cdot)$. Applying Theorem 2 to this equation we reach the
assertion
$u_t\ge 0$.

Secondly, let us prove \eqref{utposs1}. Note that
\beqst
&&D_t^\beta a(u(t,x),t,x)=\int_0^t {(t-\tau)^{-\beta}\over \Gamma(1-\beta)}{d\over d\tau}a(u(\tau,x),\tau,x)d\tau
\\
&&\quad=\int_0^t {(t-\tau)^{-\beta}\over \Gamma(1-\beta)}{d\over d\tau}[a(u(\tau,x),\tau,x)-a(u(t,x),t,x)]d\tau
\\
&&\quad =\int_0^t {\beta(t-\tau)^{-\beta-1}\over \Gamma(1-\beta)}[a(u(t,x),t,x)-a(u(\tau,x),\tau,x)]d\tau
\\
&&\quad +a(u(t,x),t,x){t^{-\beta}\over \Gamma(1-\beta)}.
\eeqst
Observing \eqref{utpossuff2} and $u_t\ge 0$ we deduce \eqref{utposs1}. Proof is complete.
\hfill $\Box$

\medskip
Verification of the relations  \eqref{ipas6},  \eqref{utpossuff2}  and interpretation of the
cone conditions \eqref{ipas4}, \eqref{ipas3a}, \eqref{utpossuff} for $z$ requires a priori information
about the bounds of $u$. Most simple way to remove $u$ from these conditions
is to assume them to hold for any $w\in \R$. For instance,
sufficient conditions for \eqref{utpossuff2} and  \eqref{ipas6} are
\beq\label{suffuld}
a(w,t,x),\, a_{w}(w,t,x),\, a_{t}(w,t,x)\ge 0
\eeq
in $\R\times (0,T)\times \Omega$ and
\beq\label{suffuld1}
\forall x\in \Omega \;\;\;\exists \varepsilon_x>0\, :\,
a(w,t,x)> 0\;
\eeq
for  $w\in\R$ and $t\in (0,\varepsilon_x)$, respectively. However, this is too restrictive, because in many practical cases the nonlinearity
function $a$ is not positive and monotone for all $w$ (see the examples below). The proved positivity principle makes possible to
restrict these conditions for positive or negative values of $w$.

Suppose that the assumptions of Theorem 2 are satisfied for $a_{ij},a_j,\omega,u$ and $a,b\in C(\R\times [0,T]\times \overline\Omega)$, $z\in
C(\overline\Omega)$.
Let us point out two cases.
\begin{description}
\item {1.}  If $u_0,g\ge 0$ and \eqref{as3} holds with $f=az+b$, then by Theorem 2 we have $u\ge 0$.
This means that for the condition  \eqref{utpossuff2}
to hold it is sufficient to assume that the relations \eqref{suffuld}
are satisfied in $(0,M)\times (0,T)\times \Omega$, where $M\in (0,\infty]$ is some {\it a priori}
known upper bound of $u$. If in addition $\forall x\in \Omega$ $\exists
\delta_x>0$ : $f(0,t,x)>0$ in $(0,\delta_x)$ then by Corollary 1 the inequality $u(t,x)>0$ is valid in  $(0,\delta_x)\times \Omega$.
Therefore, for the condition \eqref{ipas6}
to hold it is sufficient to assume that  \eqref{suffuld1} is satisfied for any $w\in (0,M)$ and $t\in (0,\varepsilon_x)$ and $\varepsilon_x$
 is taken so that
$\varepsilon_x\le \delta_x$.
\item {2.} The statements  can be reformulated for negative $u$, too. To this end, Theorem 2 and Corollary 1 have to be applied
 to the problem for $-u$.
 If $u_0,g\le 0$ and \eqref{as3} holds with $f(w,t,x)=-a(-w,t,x)z(x)-b(-w,t,x)$, then  $u\le 0$. For the condition  \eqref{utpossuff2}
to hold it is sufficient to assume that the relations \eqref{suffuld} are satisfied in $(m,0)\times (0,T)\times \Omega$, where
$m\in [-\infty,0)$ is an {\it a priori}
known lower bound of $u$.
If in addition $\forall x\in \Omega$ $\exists
\delta_x>0$ : $f(0,t,x)>0$ in $(0,\delta_x)$ then by Corollary 1 the inequality $u(t,x)<0$ is valid in  $(0,\delta_x)\times \Omega$.
For the condition \eqref{ipas6}
to hold it is sufficient to assume that  \eqref{suffuld1} is satisfied for any $w\in (m,0)$ and $t\in (0,\varepsilon_x)$ and
$\varepsilon_x\le \delta_x$.
\end{description}

Clearly, similar statements can be formulated for the conditions \eqref{ipas4}, \eqref{ipas3a}, \eqref{utpossuff}, as well.

Finally, we consider the satisfaction of the cone conditions in Theorems 3 and 4 in case of some particular
equations occurring in applications.
Let $a(w,t,x)=a(w)$ and $b(w,t,x)=b(t,x)$. This means that $za$ is an inhomogeneous reaction term and
$b$ is a source term. Let us choose the following three examples:\\
(1) linear reaction (or potential) $a(w,t,x)=w$ \cite{Jin};\\
(2) nonlinear reaction $a(w,t,x)=w(1-{w\over W})$  occurring in the Fisher equation \cite{Tar,Tur};\\
(3) nonlinear reaction $a(w,t,x)=w^2(1-{w\over W})$  in the Zeldovich equation \cite{Tar}.\\
Here $W$ is a given positive constant.

Let $u_0=0$, $g\ge 0$. Then $u\ge 0$. In context of Theorem 3 we  assume also $g_t\ge 0$ that yields $u_t\ge 0$.
Moreover, let $b(0,x)=0$ and $b(t,x)>0$ for $t\in (0,\delta)$ and some $\delta>0$. (Alternatively
we could set $b=0$, but then the additional restriction $u(t,x)>0$ for $t\in (0,\delta)$ is needed).
By means of rather elementary computations, from \eqref{ipas4}, \eqref{ipas3a}, \eqref{utpossuff2}, \eqref{utposs1}
 and  Remark 2 we deduce the following sufficient conditions in the form of direct inequalities for
 $z$ and $u$ for assumptions of Theorems 3 and 4.

The cone conditions of Theorem 4 with $\varkappa=1$ are satisfied provided\\[1ex]
$z^2\le {T^{-\beta}\over \Gamma(1-\beta)}$ in case (1);\\
$z^2\le {T^{-\beta}\over \Gamma(1-\beta)}$, $u^j\le {W\over 2}$, $j=1,2$, in case (2);\\
$z^2\le {3T^{-\beta}\over \Gamma(1-\beta)W}$, $u^j\le {2W\over 3}$, $j=1,2$, in case (3).\\[1ex]
The cone conditions of Theorem 3 hold if\\[1ex]
$z^2\le {T^{-\beta}\over \Gamma(1-\beta)}$ in case (1);\\
$z^2\le 0$, $u^j\le {W\over 2}$, $j=1,2$, in case (2);\\
$z^2\le {3T^{-\beta}\over \Gamma(1-\beta)W}$, $u^j\le {W\over 3}$, $j=1,2$, in case (3).\\[1ex]
\indent
Theorem 3 fails in case (2) for positive $z^2$. In any of the other cases, the smaller $T$,
the less restrictive the  inequalities for $z^2$ and $u^j$ are.
(For $u^j$, this follows from the homogeneous initial condition.)

\bigskip\bigskip\noindent{\bf\large Acknowledgement}

\medskip\noindent
The research was supported by Estonian Research Council Personal Research Grant PUT568 and Institutional
Research Grant IUT33-24.

\section*{Appendix}
\setcounter{section}{7}
\setcounter{equation}{0}

{\small
In this appendix we prove local  existence and global uniqueness theorems for the problem \eqref{maineq} - \eqref{bound}.

\medskip\noindent
{\bf Theorem A1}. {\it Assume that $a_{ij}$, $a_j$, $\omega$ and $p$ satisfy the assumptions of Theorem {\rm 1}. Let
 there exist a function $\widehat u\in U_p\cap C([0,T]\times\overline\Omega)$   that satisfies the
conditions\,\footnote{The space $U_p$ was defined
 in \eqref{Updefi}.}
$$\widehat u(0,x)=u_0(x),\;\; x\in\Omega,\qquad {\cal B}\widehat u(t,x)=g(t,x),\;\;(t,x)\in (0,T)\times \Omega$$ and fulfills the
relation $\Psi:=D_t^\beta (\widehat u-u_0)-A\widehat u-f(\widehat u,\cdot,\cdot)\in C([0,T];C_*(\overline\Omega))$,
where $C_*(\overline\Omega)=\{z\in C(\overline\Omega)\, :\, z|_{\partial\Omega}=0\}$
in case I and $C_*(\overline\Omega)=C(\overline\Omega)$ in case II. Concerning $f$, we assume that $f\in C(\R\times [0,T]\times\overline\Omega)$
and
\beq\label{appendix1}
\begin{array}{ll}
&\exists K,\varrho>0\; :\; |f(\widehat u+w^1,t,x)-f(\widehat u+w^2,t,x)|\le K|w^1-w^2|\;\;
\\ [1ex]
&\qquad\qquad  \forall w^1,w^2\in [-\varrho,\varrho],\; t\in [0,T],\; x\in \overline\Omega.
\end{array}
\eeq
If $T$ is smaller than a certain constant depending on $A$,  $\beta$, $\Psi$, $K$ and $\varrho$, then
the problem \eqref{maineq} - \eqref{bound} has a solution in the space $U_p\cap C([0,T]\times\overline\Omega)$.
}
\\[1ex]
{\it Proof}. Firstly, let us show that the following  inequality is valid:
\beq\label{appendix2}
\begin{array}{ll}
&\|J_t^\beta y\|_{C^\alpha([0,T];X)}\le {1\over\beta\Gamma(\beta)}(T^\beta+2T^{\beta-\alpha})\|y\|_{C([0,T];X)}
\end{array}
\eeq
for any $y\in C([0,T];X)$,
where $X$ is a Banach space and $\alpha\in (0,\beta)$. Indeed,
\beqst
\begin{array}{ll}
&\Gamma(\beta)\|J_t^\beta y\|_{C^\alpha([0,T];X)}=\sup\limits_{0<t<T}\bigl\|\int_0^t (t-\tau)^{\beta-1} y(\tau)d\tau \bigl\|
\\[1ex]
&\quad + \sup\limits_{0<s<t<T}{1\over (t-s)^\alpha} \bigl\|\int_0^t (t-\tau)^{\beta-1} y(\tau)d\tau - \int_0^s (s-\tau)^{\beta-1} y(\tau)d\tau\bigl\|
\\[2ex]
&\le (I_1+I_2+I_3)\|y\|_{C([0,T];X)},
\end{array}
\eeqst
where
\beqst
\begin{array}{ll}
&I_1=\sup\limits_{0<t<T}\int_0^t (t-\tau)^{\beta-1}d\tau=\sup\limits_{0<t<T}\int_0^t \tau^{\beta-1}d\tau={T^\beta\over \beta},
\\[2ex]
&I_2=\sup\limits_{0<s<t<T}{1\over (t-s)^\alpha}\int_0^s \bigl|(t-\tau)^{\beta-1}-(s-\tau)^{\beta-1}\bigl| d\tau
=\sup\limits_{0<s<t<T}{1\over (t-s)^\alpha}
\\  [2ex]
&\times \bigl(\int_0^s \tau^{\beta-1}d\tau - \int_{t-s}^t \tau^{\beta-1}d\tau\bigl)\le
\sup\limits_{0<s<t<T}{1\over (t-s)^\alpha} \int_0^{t-s}\tau^{\beta-1}d\tau={T^{\beta-\alpha}\over \beta},
\\ [2ex]
&I_3=\sup\limits_{0<s<t<T}{1\over (t-s)^\alpha}\int_s^t(t-\tau)^{\beta-1}d\tau
=\sup\limits_{0<s<t<T}{1\over (t-s)^\alpha}\int_0^{t-s}\tau^{\beta-1}d\tau={T^{\beta-\alpha}\over \beta}.
\end{array}
\eeqst
This implies \eqref{appendix2}.

Making use of the change of variables $u=\widehat u+v$,  the problem \eqref{maineq} - \eqref{bound} for $u\in U_p\cap C([0,T]\times\overline\Omega)$
is reduced to the following problem  for $v\in U_p\cap C([0,T]\times\overline\Omega)$:
\beq\label{appendix3}
\begin{array}{ll}
&D_t^\beta v(t,x) = A(x)v(t,x) +
F(v(t,x),t,x)\, ,\quad t\in (0,T),\; x\in \Omega\, ,
\\ [1ex]
&v(0,x)=0\, ,\quad x\in\Omega, \quad
{\cal B}v(t,x)=0\, ,\;\; (t,x)\in (0,T)\times {\partial\Omega},
\end{array}
\eeq
where $F(v,t,x)=f(\widehat u+v,t,x)-f(\widehat u,t,x)-\Psi(t,x)$.
Due to Lemma 1, \eqref{appendix3} is in $U_p\cap C([0,T]\times\overline\Omega)$ equivalent to the problem
\beq\label{appendix4}
\begin{array}{ll}
&v(t,x) =A(x) J_t^\beta v(t,x) +
J_t^\beta F(v(t,x),t,x)\, ,\quad t\in (0,T),\; x\in \Omega\, ,
\\ [1ex]
&v(0,x)=0\, ,\quad x\in\Omega, \quad
{\cal B}v(t,x)=0\, ,\;\; (t,x)\in (0,T)\times {\partial\Omega}.
\end{array}
\eeq

On the other hand, there exists $\xi\in\R$ such that the operator ${\cal A}=A+\xi$ satisfies the assumptions of Lemma 4 in the space
$C_*(\overline\Omega)$ with the domain $D({\cal A})=\{z\, :\, z\in W_q^2(\Omega)\, \forall q\in (1,\infty), {\cal A}z\in C_*(\overline\Omega),\,
 {\cal B}z|_{\partial\Omega}=0\}$
(see \cite{Luna}, Sect. 3.1.5). Let us consider the following operator equation:
\beq\label{appendix5}
v={\cal F}v \qquad \mbox{with\,\,\, ${\cal F}={\cal Q}_\alpha J_t^\beta [F(v,\cdot,\cdot)-\xi v]$}
\eeq
in the ball $B_\varrho=\{v\, :\, v\in C([0,T];C_*(\overline\Omega)), \, \|v\|_{C([0,T];C_*(\overline\Omega))}\le\varrho\}$, where
$\alpha=\beta/2$ and ${\cal Q}_\alpha$ is the operator defined in Lemma 4.
Using Lemma 4, the definition of $F$,  \eqref{appendix1} and \eqref{appendix2} we obtain for any $v\in B_\varrho$
the estimate
\beqst\label{appendix6}
\begin{array}{ll}
&\|{\cal F}v\|_{C([0,T];C_*(\overline\Omega))}\le \|{\cal F}v\|_{C^{\alpha}([0,T];C_*(\overline\Omega))}
\le {\|{\cal Q}_{\alpha}\|\over\beta\Gamma(\beta)}(T^\beta+2T^{\beta/2})
\\[1ex]
&\quad\times \bigl[(K+|\xi|)\|v\|_{C([0,T];C_*(\overline\Omega))}+
\|\Psi\|_{C([0,T];C_*(\overline\Omega))}\bigl]
\\[1ex]
&\quad\le   {\|{\cal Q}_\alpha\|\over\beta\Gamma(\beta)}(T^\beta+2T^{\beta/2})
 \bigl[(K+|\xi|)\varrho+
\|\Psi\|_{C([0,T];C_*(\overline\Omega))}\bigl].
\end{array}
\eeqst
From this estimate we see that if $T$ is smaller than a certain constant $T_1$  depending on $\|{\cal Q}_\alpha\|$, $\beta$, $K$, $\xi$, $\varrho$ and
$\Psi$, then the inequality $\|{\cal F}v\|_{C([0,T];C_*(\overline\Omega))}\le\varrho$ is valid. This means that for  $T<T_1$ the operator
${\cal F}$ maps $B_\varrho$ into $B_\varrho$. Similarly, for any
$v^1,v^2\in B_\varrho$ we deduce the estimate
$$
\begin{array}{ll}
&\|{\cal F}v^1-{\cal F}v^2\|_{C([0,T];C_*(\overline\Omega))}\le \kappa \|v^1-v^2\|_{C([0,T];C_*(\overline\Omega))} \end{array}
$$
with $\kappa={\|{\cal Q}_\alpha\|\over\beta\Gamma(\beta)}(T^\beta+2T^{\beta/2})
(K+|\xi|)$.
There exists $T_2$ depending on $\|{\cal Q}_\alpha\|$,  $\beta$, $K$, $\xi$ such that if $T<T_2$ then
$\kappa<1$. We have shown that ${\cal F}$ is a contraction in $B_\varrho$ provided $T<T_3=\min\{T_1;T_2\}$. Thus,
\eqref{appendix5} has a unique solution in $B_\varrho$ for  $T<T_3$.

Further, for this solution we have the relation $F(v,\cdot,\cdot)-\xi v\in C([0,T];C_*(\overline\Omega))$. Therefore, in view of
\eqref{appendix2}, it holds $\varphi=J_t^\beta [F(v,\cdot,\cdot)-\xi v]\in C^\alpha([0,T];C_*(\overline\Omega))$. By virtue
of Lemma 4, the operator equation \eqref{appendix5} is equivalent to \eqref{appendix4} and its solution $v$ belongs to
$C^\alpha([0,T];C_*(\overline\Omega))$. On the other hand,
since $[F(v,\cdot,\cdot)-\xi v]\in L_p((0,T);L_p(\Omega))$, due to Lemma 1 we get $J_t^\beta [F(v,\cdot,\cdot)-\xi v]\in
H_p^\beta ([0,T];L_p(\Omega))$ and Lemma 2 implies that the solution $v$ of \eqref{appendix4} belongs to $U_p$. Summing up,
$v\in U_p\cap C([0,T];C_*(\overline\Omega))$ and in view of the assumptions imposed on $\widehat u$ we obtain
$u=\widehat u+v\in U_p\cap C([0,T]\times \overline\Omega)$. The proof is complete. \hfill $\Box$
\\[1ex]
{\bf Remark A1}. Let us briefly discuss simple possibilities to construct functions $\widehat u$ satisfying the conditions of
Theorem A1. Firstly, consider the case II. Assume
$u_0\in C^{2+l}(\overline\Omega)$,
$g\in C^{1+l,{1+l\over 2}}([0,T]\times\partial\Omega)$ for some $l\in (0,1)$, where
$C^{\gamma,{\gamma\over 2}}$ is the anisotropic H\"older space (see \cite{Ladyz}), and ${\cal B}u_0(x)=g(0,x)$ for $x\in \partial\Omega$.
Define  $\widehat u$ as the solution of
$\widehat u_t-A\widehat u=0$ subject to \eqref{ini} and \eqref{bound}. Then
$\widehat u$, $\widehat u_t$ and $A\widehat u$ are continuous in $[0,T]\times \overline\Omega$ \cite{Ladyz}, hence $D_t^\beta (\widehat u-u_0)=
J_t^{1-\beta}\widehat u_t$ and $f(\widehat u,\cdot,\cdot)$ are also continuous. We immediately obtain $\Psi\in C([0,T];C_*(\overline\Omega))$.
In case I, let $u_0\in C^{2+l}(\overline\Omega)$,
$g\in C^{2+l,1+{l\over 2}}([0,T]\times\partial\Omega)$ for some $l\in (0,1)$ and $u_0(x)=g(0,x)$ for $x\in \partial\Omega$.
Moreover, assume $f(g,\cdot,\cdot)\in C^{l,{l\over 2}}([0,T]\times\partial\Omega)$. Then the function
$\zeta=g_t -D_t^\beta (g-u_0)+f(g,\cdot,\cdot)=g_t - J_t^\beta g_t+f(g,\cdot,\cdot)$ belongs to $C^{l,{l\over 2}}([0,T]\times\partial\Omega)$.
Let us continue $\zeta$ to a function $C^{l,{l\over 2}}([0,T]\times\overline\Omega)$ and assume additionally
$Au_0(x)+f(g(0,x),0,x)=0$ for $x\in\partial\Omega$.
Define  $\widehat u$ as the solution of
$\widehat u_t-A\widehat u=\zeta$ subject to \eqref{ini} and \eqref{bound}. Again,
$\widehat u$, $\widehat u_t$ and $A\widehat u$ are continuous, hence $\Psi$ is continuous. Due to the choice of $\zeta$ we
have $\Psi(t,x)=0$ for $x\in\partial\Omega$. Thus, $\Psi\in C([0,T];C_*(\overline\Omega))$.
\\[1ex]
{\bf Theorem A2}. {\it Assume that the assumptions of Theorem {\rm 2} are satisfied for $a_{ij}$, $a_j$, $\omega$ and \eqref{as2}
hold.
Suppose that \eqref{maineq} - \eqref{bound} has two solutions $u^1$ and $u^2$ satisfying the relations
$u^l \in C^{\beta'}([0,T];C(\overline\Omega))$, $l=1,2,$ for some $\beta'>\beta$ and
 $u^l_{x_j},u^l_{x_i,x_j}\in C((0,T]\times \overline\Omega)$, $l=1,2$.
Moreover, let the inequality $|f(u^1(t,x),t,x)-f(u^2(t,x),t,x)|\le M|u^1(t,x)-u^2(t,x)|$ be valid for
$t\in [0,T]$, $x\in \overline \Omega$ with some $M>0$.
Then $u^1=u^2$.
}
\\[1ex]
{\it Proof}. The difference $u=u^1-u^2$ solves the problem
\beq\label{appendix10}
\begin{array}{ll}
&D_t^\beta u(t,x) = A(x)u(t,x) +
G(u(t,x),t,x)\, ,\quad t\in (0,T),\; x\in \Omega\, ,
\\ [1ex]
&u(0,x)=0\, ,\quad x\in\Omega, \quad
{\cal B}u(t,x)=0\, ,\;\; (t,x)\in (0,T)\times {\partial\Omega},
\end{array}
\eeq
with linear with respect to $w$ function
 $G(w,t,x)={f(u^1(t,x),t,x)-f(u^2(t,x),t,x)\over u^1(t,x)-u^2(t,x)} w$. Applying Theorem 2 to problems of $u$ and $-u$ we obtain
$u\ge 0$ and $u\le 0$. This implies the assertion of the theorem. \hfill $\Box$
}

\end{document}